\documentclass[]{amsart}
\usepackage{color}

 \usepackage[latin1]{inputenc}
 \usepackage[T1]{fontenc}
 \usepackage[normalem]{ulem}
\usepackage{verbatim}
 \usepackage{graphicx}
\usepackage[latin1]{inputenc}
\usepackage{amsmath}
\usepackage{amssymb}
\usepackage{amsfonts}
\usepackage{amsthm}
\usepackage{bbm}

\newcommand{\C}{\mathbb{C}}
\newcommand{\N}{\mathbb{N}}
\newcommand{\R}{\mathbb{R}}

\renewcommand{\Im}{{\rm Im}}

\newtheorem{theorem}{Theorem}[section]

\theoremstyle{remark}
\newtheorem{remark}{Remark}[section]
\theoremstyle{definition}

\begin{document}

\title[Davey-Stewartson]{A numerical 
approach to Blow-up issues for Davey-Stewartson II type systems}
\author{Christian Klein}
\address{Institut de Math\'ematiques de Bourgogne\\
                Universit\'e de Bourgogne, 9 avenue Alain Savary, 21078 Dijon
                Cedex, France\\
    E-mail Christian.Klein@u-bourgogne.fr}

\author{Jean-Claude Saut}
\address{Laboratoire de Math\' ematiques, UMR 8628\\
Universit\' e Paris-Sud et CNRS\\ 91405 Orsay, France\\ E-mail jean-claude.saut@math.u-psud.fr}
\date{June 6th, 2014}
\maketitle
\large

\begin{abstract}
We provide a numerical study of various issues pertaining to the 
dynamics of the Davey-Stewartson systems of the DS II type. In particular we investigate whether or not the properties (blow-up, radiation,...) displayed by the focusing and defocusing DS II integrable systems persist in the non integrable case.
\end{abstract}

\vspace{0.3cm}
\begin{center}
{\it To Gustavo Ponce with friendship and admiration}
\end{center}

\section{Introduction}
This paper is concerned with blow-up issues and the long-time behavior of solutions to 
 Davey-Stewartson (DS) II type systems,
\begin{equation}
    \label{DSII}
\begin{array}{ccc}
i\epsilon 
\partial_{t}\psi+\epsilon^{2}\partial_{xx}\psi-\epsilon^{2}\partial_{yy}\psi+2\rho\left(\beta\phi+\left|\psi\right|^{2}\right)\psi & = & 0,\quad \psi :  \R^2\times \R \to \C,
\\
\partial_{xx}\Phi+\partial_{yy}\Phi+2\partial_{xx}\left|\psi\right|^{2} & = & 0, \quad \Phi :  \R^2\times \R \to \R,
\end{array}
\end{equation}
where  $\rho$ takes the values $\pm1$, and where
$\epsilon$ is  a dispersion
parameter which depending on the circumstances may be small or of order one.  Since $\epsilon$ has the 
same role as the $\hbar$ in the Schr\"odinger equation, the limit 
$\epsilon\to0$ is also called the semiclassical limit in this 
context. It can be introduced in the DS system for $\epsilon=1$ via a 
transformation $x\to x/\epsilon$, $y\to y/\epsilon$, $t\to t/\epsilon$. Since we are in particular
interested in the study of the long-time behavior of solutions, this 
is equivalent to the small $\epsilon$ behavior in $t$. At the same 
time, small $\epsilon$ implies the study of solutions on large scales 
of $x,y$. It turns out that this limit $\epsilon\to0$ can be 
conveniently handled numerically. The alternative would be to study 
the system (\ref{DSII}) for long times on larger domains, which is of 
course equivalent numerically to the case of small $\epsilon$. 

The general Davey-Stewartson systems \eqref{DSgen} are a simplification of the Benney-Roskes, Zakharov-Rubenchik systems (\cite {BR, ZR}) who in turn 
are "universal" models for the description of interaction of short 
and long waves. They were first derived in the context of water waves (\cite {DS, DR, AS}) in the 
so-called modulational (Schr\"{o}dinger) regime, and then $\Phi$ is the  mean flow induced by the interaction of oscillating modes (see \cite{La} for more details and references and for a rigorous derivation of the Davey-Stewartson systems  for water waves).

They were also rigorously proven in \cite{Co,CoLa} to provide a good approximate solution to general quadratic hyperbolic systems using diffractive geometric optics. In fact they have been formally derived in many concrete physical contexts, ferromagnetism \cite{Le}, plasma physics \cite{MRZ}, nonlinear optics \cite{NM}.

The Davey-Stewartson systems can also be viewed as  the two-dimensional version  of the {\it Zakharov-Schulman} 
systems (see \cite{ZS, ZS2, GS3}) which read in dimension $d=2, 3:$

\begin{equation}\label{ZS}
\begin{array}{ccc}
i\partial_{t}\psi +L_1 \psi+\psi \Phi=0,\\
L_2\Phi=L_3|\psi|^2,
\end{array}
\end{equation}

where

$$L_n=\sum_{j,k=1}^d C_{jk}^n\frac{\partial^2}{\partial x_j\partial x_k}, \; n=1,2,3.$$

are second order differential operators with constant coefficients, 
the matrices $\left(C_{j,k}^n\right )_{1\leq j\leq k\leq d}$ being real and symmetric (but not necessarily positive or negative). When $d=2$,  the Zakharov-Schulman systems reduce in fact to Davey-Stewartson systems. It is worth noticing that Schulman \cite{Sc} has proven that among the systems obtained in the $d=2$ and $d=3$ cases, the only ones which are integrable are those previously known in the two-dimensional cases, that is the DS I and DS II systems (see below). In particular, none of them are integrable in the three-dimensional case.

It will be convenient later on to write \eqref{ZS} in a slightly different form, introducing $d$ real valued functions $\phi_1,...,\phi_d$ satisfying

\begin{equation}\label {phitwo}
L_2\phi_j=\frac{\partial |\psi|^2}{\partial x_j},\quad j=1,...,d.
\end{equation}

We set (for $\phi=(\phi_1,...,\phi_d)\in \R^d)$

$$\mathcal L_3\phi =\sum_{j,k=1}^d C_{jk}^3 \frac{\partial \phi_j}{\partial x_k},$$

which allows to rewrite \eqref{ZS}, assuming that the matrix 
$C^2_{ij}$ is invertible, as

\begin{equation}\label{ZSbis}
i\partial_{t}\psi +L_1 \psi+(\mathcal L_3 \phi )\psi=0.
\end{equation}

\vspace{0.3cm}
The Davey-Stewartson systems have in fact the general form (when using  the variable $\phi$ such that $\phi_x=\Phi $ instead of $\Phi$ in \eqref{DSII}) , where $a, b, c, \nu_1, \nu_2$ are real parameters depending on the physical context

\begin{equation}\label{DSgen}
\begin{array}{ccc}
i\partial_{t}\psi +a\partial_x^2\psi+b\partial_y^2 \psi=(\nu_1|\psi|^2+\nu_2\partial_x\phi)\psi,\\
c\partial_x^2\phi+\partial_y^2\phi=-\delta\partial_x|\psi|^2,
\end{array}
\end{equation}
where one can assume (up to a change of unknown) $b>0$ and $\delta>0.$

Using the terminology of \cite{GS}, one says that \eqref{DSgen} is

\begin{center}

elliptic-elliptic if \quad ($\text{sgn}\; a, \text{sgn}\; c)=(+1,+1),$

hyperbolic-elliptic  if \quad ($\text{sgn}\; a, \text{sgn}\; c)=(-1,+1),$

elliptic-hyperbolic  if \quad ($\text{sgn}\; a, \text{sgn}\; c)=(+1,-1),$

hyperbolic-hyperbolic  if \quad ($\text{sgn}\; a, \text{sgn}\; c)=(-1,-1).$

\end{center}

The so-called {\it DS I and DS II systems} are integrable, but very 
particular cases of respectively elliptic-hyperbolic and 
hyperbolic-elliptic Davey-Stewartson systems. In fact they correspond 
to a very special choice of the coefficients in \eqref{DSII} or 
\eqref{DSgen} that  have  limited  physical relevance (in the context of water waves they occur in the shallow water limit, see \cite{AS}). One may hope however that they will give some insights into the dynamics of the corresponding non integrable systems for which standard PDE tools provide only local well-posedness results (global for small data) without any qualitative information.

The hyperbolic-hyperbolic case does not seem to occur in a physical 
situation. The elliptic-elliptic DS systems may possess blow-up 
solutions by focusing (see \cite{GS}), similarly to the focusing  
cubic NLS\footnote{Note however that a precise analysis of the 
blow-up as in \cite{MeRa} for the critical focusing NLS seems to be 
missing.}. Numerical simulations can be found in \cite{PSSW, BMS}. 
The elliptic-hyperbolic case offers the more challenging mathematical 
and numerical problems. The local Cauchy problem has been first 
studied by Linares and Ponce in \cite{LiPo2} (see  \cite{NH1,H} for the best known results 
so far). The global existence of small solutions with the large 
time   asymptotic) was proven in \cite{NH2}.  We refer to   \cite{BB,BMS,MFP} for numerical simulations. 

A precise asymptotic of small solutions to the {\it integrable} DS I system is given in \cite{Ki}.

The DS II type systems \eqref{DSII} are thus the hyperbolic-elliptic ones. They are in particular  relevant for surface gravity waves without surface tension (see \cite{La} for a rigorous justification).

The DS II type systems for water waves  reduce in the infinite depth limit  to the "hyperbolic" nonlinear Schr\"{o}dinger equation

\begin{equation}\label{HNLS}
i\partial_{t}\psi+\partial_{xx}\psi-\partial_{yy}\psi+|\psi|^2\psi=0
\end{equation}
which was derived in \cite{Za} as a model of gravity waves in the modulation regime and in infinite depth.

\vspace{0.3cm}
The present paper will be uniquely devoted to the DS II type systems \eqref{DSII}. The parameter $\beta\in \mathbb{R}$ in \eqref{DSII} determines the 
contribution of the mean field $\Phi$ to the nonlinearity in the 
nonlinear Schr\"odinger (NLS) 
equation. For $\beta=0$ one obtains the hyperbolic NLS equation, for 
$\beta=1$ a completely integrable system.  There is a focusing ($\rho=-1$) and a 
defocusing ($\rho=1$) version of DS II type systems except for $\beta=0$. In the 
case of the hyperbolic NLS, there is always one focusing and one 
defocusing direction. A change of sign of $\rho$ in this case just 
interchanges (up to complex conjugation) the role of the spatial variables.

The DS II system can be viewed as a nonlocal cubic nonlinear Schr\"{o}dinger equation. Actually one can solve $\Phi$ as

$$\Phi=2\lbrack(-\Delta)^{-1}\partial_{xx}\rbrack|\psi|^2,$$

where $(-\Delta)^{-1}\partial_{xx}=R_1^2$ is a zero order operator with Fourier symbol $-\frac{\xi_1^2}{|\xi|^2}$ and is thus bounded in all $L^p(\R^2)$ spaces, $1<p<\infty$ and all  Sobolev spaces $H^s(\R^2),$ allowing to write \eqref{DSII} as

\begin{equation}\label{DDSIIeq}
i\epsilon 
\partial_{t}\psi+\epsilon^{2}\partial_{xx}\psi-\epsilon^{2}\partial_{yy}\psi+2\rho\left(2\beta R_1^2 (|\psi|^2)+\left|\psi\right|^{2}\right)\psi =  0.
\end{equation}

One easily finds that  (\ref{DSII}) has two formal conservation laws, the $L^2$ norm

$$ \int_{\R^2}|\psi(x,y,t)|^2dx dy= \int_{\R^2}|\psi(x,y,0)|^2dx dy$$

and the energy (Hamiltonian)
\begin{eqnarray}
    E(\psi(t))&=&\int_{\R^2}\left[\epsilon^2 
|\partial_{x}\psi|^2-\epsilon^2|\partial_{y}\psi|^2-\rho (|\psi|^2+ 
\beta\Phi)|\psi|^2)\right] dxdy\nonumber\\
&=&E(\psi(0)).
    \label{energy}
\end{eqnarray}

\vspace{0.3cm}
One can also write  (\ref{DSII}) as

\begin{equation}\label{DSIIbis}
   i \epsilon\partial_{t}\psi+\epsilon^{2}\partial_{xx}\psi-\epsilon^{2}\partial_{yy}\psi+2\rho\Delta^{-1}\lbrack
    \left(\partial_{yy}+(1-2\beta)\partial_{xx}\right)\left|\psi\right|^{2}\rbrack\psi  =  0,
    \end{equation}
    
 which involves the order zero nonlocal operator

$$\Delta^{-1}\lbrack
    \left(\partial_{yy}+(1-2\beta)\partial_{xx}\right)\rbrack.$$
    

We will in the following computations always study initial data in the Schwartz 
space of rapidly decreasing functions. Numerically these will be 
treated as essentially periodic which allows the use of Fourier 
spectral methods. Thus we solve the DS II system in the form 
(\ref{DSIIbis}) in Fourier space. 

Note that the integrable case $\beta=1$ is 
distinguished by the fact that the same hyperbolic operator appears 
in the linear and in the nonlinear part. In this case 
the equation is invariant under the transformation $x\to y$ and $\psi 
\to \bar{\psi}$ and \eqref{DSIIbis} can be written in a "symmetric" form as

\begin{equation}\label{box}
i \epsilon\partial_{t}\psi+\epsilon^{2}\Box\psi-2\rho\lbrack(\Delta^{-1}\Box) |\psi|^2\rbrack\psi=0,
\end {equation}

where $\Box =\partial_{xx}-\partial_{yy}.$

This extra symmetry in the integrable case could be responsible for 
the existence of localized lump solutions and to the blow-up phenomena that do not not  seem persist in the non integrable case.

\vspace{0.3cm}
The paper will be organized as follows. We first recall some rigorous known results on the Cauchy problem for the DS II type systems. It turns out that Inverse Scattering techniques 
allow to provide very precise information on the dynamics of the DS 
II systems, both in the integrable focusing and defocusing cases. One aim of the 
present paper is to investigate numerically to what extent those 
dynamics are generic in the sense that they hold also in the non 
integrable cases. This issue motivates the numerical simulations of 
Section 3 which suggest that  the blow-up  which occurs in  the 
focusing DS II system  does not persist in the (focusing) non integrable case, while the purely dispersive regime of the defocusing DS II system seems to persist in the (defocusing) non integrable case.

\subsection*{Notations}

The following notations will be used throughout this article. The partial derivative with be denoted by $\partial_x \phi,..$. For any $s\in \R,$ $D^s=(-\Delta)^{\frac{s}2}$ and $J^s=(I-\Delta)^{\frac{s}2}$ denote the Riesz and Bessel potentials of order $-s$, respectively. 

The Fourier transform of a function $f$ is denoted by $\hat{f} $ or 
$\mathcal{F}(f).$ For $1 \le p \le \infty$, $L^p(\mathbb R)$ is the 
usual Lebesgue space with the norm $|\cdot |_p$, and for $s \in 
\mathbb R$, the Sobolev spaces $H^s(\mathbb R^2)$ are defined via the usual norm $\|\phi \|_s= |J^s \phi|_2$. 

 $\mathcal S(\R^2)$ will denote the Schwartz spaces of smooth rapidly decaying functions, and $\mathcal S'(\R^2)$ the space of tempered distributions.

\section{Summary of known theoretical results}

\subsection{The Cauchy problem}

As previously noticed the DS II type systems can be reduced to   a nonlocal cubic (hyperbolic) nonlinear Schr\"{o}dinger equation (HNLS), and the local Cauchy problem
theory for the cubic NLS extends  to the DS II systems. Namely one has the following local well-posedness result (\cite{GS}):

\begin{theorem}\label{GS}
(i) Let $\psi_0\in L^2(\R^2).$ There exists a unique maximal solution 
$(\psi,\phi)$ of \eqref{DSII} on $\lbrack 0,T^*)$, $T^*>0,$ such that
$$\psi\in C(\lbrack 0,T^*), L^2(\R^2))\cap L^4((0,t)\times \R^2),$$
$$\Phi\in L^2((0,t)\times \R^2),$$
$$\psi(0)=\psi_2,\quad |\psi(t)|_2=|\psi_0|_2,\quad 0\leq t<T^*.$$

(ii) If $|\psi_0|_2$ is sufficiently small, then $T^*=+\infty,$ the solution is global.

(iii) If $\psi\in H^1(\R^2),$ the previous solution satisfies
$$\psi\in C(\lbrack 0,T^*); H^1(\R^2))\cap C^1(\lbrack 0,T^*); 
H^{-1}(\R^2)),$$
$$\nabla \psi\in L^4((0,t)\times \R^2), \quad \Phi\in C(\lbrack 0,T^*);L^p(\R^2)),$$
$$\nabla \Phi \in L^4(0,t;L^q(\R^2)),$$
\vspace{0.3cm}
for every $t\in \lbrack 0,T^*),\; p\in \lbrack 2,\infty)$\; and\; $q\in \lbrack 2,4\rbrack.$
\vspace{0.3cm}\\
Moreover, $E(t)=E(0), 0\leq t<T^*.$

\end{theorem}

\begin{remark}
 (i) Results in Theorem  \ref{GS} do not distinguish between a {\it 
 focusing} and a {\it defocusing case}. In fact they rely only on the 
 dispersive (Strichartz) estimates for the linear group 
 $e^{it(\partial_{xx}-\partial_{yy})}$ which are the same as those of 
 the standard Schr\"{o}dinger group (see \cite{GS4}). Since 
 \eqref{DSIIbis} is $L^2$ critical, the maximal existence time $T^*$ does not depend only on the $L^2$ norm of the initial data $\psi_0$, but on $\psi_0$ in  a more complicated way (see \cite{CW}). This explains why the conservation of the $L^2$ norm of the solution does not imply global well-posedness.
 
 (ii) The global existence result in part (ii) of Theorem \ref{GS} does not provide any large time asymptotic. 
 
(iii)  In the $H^1$ context, the conservation of energy does not yield any $H^1$ bound since its quadratic part is not positive definite.

(iv) All the previous considerations apply as well to  the hyperbolic NLS  \eqref{HNLS} for which global existence  with scattering 
is conjectured (see the numerical simulations below).

\end{remark}





\vspace{0.3cm}
More can be said in the integrable case ($\beta =1),$ where the addition of   a nonlocal cubic term to the hyperbolic NLS equation may have dramatic effects in the focusing case. Moreover the integrability provides in the defocusing case results that would be difficult to obtain by pure PDE techniques. In particular, Sung (\cite{Su1, Su2, Su3, Su4}) has proven the following

\begin{theorem}\label{DSSung}
Let $\psi_0\in \mathcal S(\R^2)$\footnote{This condition can 
obviously be weakened.}. Then \eqref {DSII} possesses a unique global solution $\psi$ such that the mapping $t\mapsto \psi(\cdot,t)$ belongs to $C^\infty(\R, \mathcal S(\R^2))$ in the two cases:

(i) Defocusing.

(ii) Focusing and $|\widehat{\psi_0}|_1|\widehat{\psi_0}|_\infty<C,$ where $C$ is an explicit constant.

Moreover, there exists $c_{\psi_0}>0$ such that

$$|\psi(x,t)|\leq \frac{c_{\psi_0}}{|t|}, \quad (x,t)\in \R^2\times \R^*.$$
\end{theorem}

We recall that such a result is unknown for the general {\it non integrable } DS-II systems, and also for the hyperbolic cubic NLS.

\begin{remark}
1. Sung obtains in fact the global well-posedness (without the decay rate) in the defocusing case under the assumption that $\hat \psi_0\in L^1(\R^2)\cap L^\infty(\R^2)$ and $\psi_0\in L^p(\R^2)$ for some $p\in [1,2),$ see \cite{Su4}.

2. Recently, Perry \cite{Pe} has given a more precise asymptotic behavior in the defocusing case for initial data in $H^{1,1}(\R^2)=\lbrace f\in L^2(\R^2) \;\text{such that}\; \nabla f, (1+|\cdot|)f \in L^2(\R^2)\rbrace,$ proving that the solution obeys the asymptotic behavior in the $L^\infty(\R^2)$ norm : 

$$\psi(x,t)=u(\cdot,t)+o(t^{-1}),$$

where $u$ is the solution of the linearized problem. We should emphasize again that this kind of result is out of reach of the present PDE techniques. As previously said, one of the aims of the present paper is to give numerical evidence that this behavior persists in the non integrable defocusing case.

\end{remark}

On the other hand, the integrable focusing  DS II system possesses a family of localized solitary waves (\cite{APP, AC}), the lumps

\begin{equation}\label{lump}
 \psi(x,y,t) = 2c \frac{\exp \left( -2i(\xi x - \eta y + 
2(\xi^{2}-\eta^{2} )t)\right)}{|x + 4\xi t + i(y + 4\eta t) +
 z_{0}|^2+|c|^2}
\end{equation}
 where $(c,z_{0})\in \mathbb{C}^2$ and $(\xi,\eta)\in\mathbb{R}^2$ 
 are constants. The lump moves with constant 
 velocity $(-4\xi, -4\eta)$ and decays as $(x^2+y^2)^{-1}$ for 
 $x,y\to\infty$.

  Again in  the  focusing integrable case, Ozawa \cite{Oz} has constructed an explicit  solution of the Cauchy problem  whose $L^2$ norm blows up in finite time $T^*$. In fact the mass density $|\psi(.,t)|^2$ of the solution converges as $t\to T^*$ to a Dirac measure with total  mass $|\psi(.,t)|^2_2=|\psi_0|^2_2$  (a weak form of the conservation of the $L^2$ norm). Every regularity breaks down at  the blow-up point but the solution persists after the blow-up time and disperses in the sup norm when  $t\to \infty $ as $t^{-2}.$ The construction of the blow-up solution is obtained by  applying the pseudo-conformal invariance law which holds for Davey-Stewartson systems (see \cite{GS1}) to  a  lump solitary wave (see above). More precisely,

\begin{theorem}[Ozawa] 
Let   $ab<0$ and $T=-a/b$. 
Denote by $u (x,y,t)$ the function defined by
\begin{equation}
 \psi (x,y,t) = \exp \left( i \frac{b}{4(a+bt)} (x^2 - y^2)
\right) \frac{v(X,Y)}{a+bt}
\label{soluoz}
\end{equation}
 where 
\begin{equation}
 v(X,Y) = \frac{2}{1+X^2+Y^2}, \,\, X=\frac{x}{a+bt},
\,\, Y=\frac{y}{a+bt}
\end{equation}
Then, $\psi$ is a solution of (\ref{DSII}) with 
\begin{equation}
 \|\psi(x,y,t)\|_2 = \|v(X,Y)\|_2 = 2\sqrt{\pi}
\end{equation}
and
\begin{equation}
 |\psi(t)|^2 \rightarrow 2\pi \delta 
 \,\, \text{in}\quad \mathcal S'\quad \mbox{when} \,\, t \rightarrow T.
\end{equation}
where $\delta$ is the Dirac measure.
\end{theorem}



Note that this construction is reminiscent of a similar one for the $L^2$ critical focusing NLS in $\R^n$

\begin{equation}
\label{criNLS}
i\partial_{t}\psi+\Delta \psi+|\psi|^{n/2}\psi =0,
\end{equation}

constructed via a pseudo-conformal transformation  applied to  a ground state solution of NLS.

It is well known (\cite{MeRa} and the references therein) that this blow-up is non generic for \eqref{criNLS} (it does not give the generic blow-up rate).  

However the blow-up shown by Ozawa is very different from the NLS 
one. First, it is an $L^2$ blow-up and not an $H^1$ one ($u(.,.,0)\notin H^1(\R^2)$). 
Lastly, the solution extends beyond the blow-up time and then scatters as $t\to\pm \infty.$ More precisely, it is shown in \cite{Oz} that there  exists a unique (explicit) $f\in L^2(\R^2)$ such that 

\begin{equation}\label{scatt}
|u(t)-U(t)f|_2\to 0\quad \text{as}\quad t\to \pm \infty,
\end{equation}

where $U(t)$ denotes the unitary group $U(t)=\text{exp}(it(\partial_x^2-\partial_y^2)).$

Recall that no blow-up occurs for DS II type equation with initial data small enough in $L^2$but a precise $L^2$ bound is not known.
\footnote{Recall that for the focusing cubic equation in 2D, the criterion 
is $|\psi_0|_2<|Q|_2,$ where $Q$ is the ground state solution of the focusing $L^2$ critical  NLS.}. There is such a criterion in Theorem \ref{DSSung} (this is not an $L^2$ one) and one can ask whether or not it is optimal in Sung functional setting.

We do not know of a proof of a blow-up for the focusing DS II 
equation occurring for initial data different from the above construction (say a Gaussian with sufficiently large mass).  Our numerical simulations suggest that blow-up actually  may occur in this situation which is not theoretically well understood.

The "Ozawa blow-up"   is carefully studied numerically in 
\cite{KMR}   illustrating in particular that the Ozawa solution is 
unstable. The same "structural instability" is shown as in \cite{MFP} 
for the DS II lump (see \cite{Ki}). On the other hand, our numerical 
simulations below suggest that neither an  Ozawa type blow-up nor 
another type of generic blow-up  persist in the non integrable case. 

Since Ozawa's blow-up is obtained by using both the existence of a lump solution and of a pseudo-conformal law, one may first ask whether such a transformation exist for general DS type systems. In fact it was established in \cite{GS3} that this is the indeed the case for the more general Zakharov-Schulman systems \eqref{ZS}. For instance, assuming that the matrix $C^1$ is invertible and denoting by $q$ the following quadratic form on $\R^2$:

$$q(x)=\sum_{j,k=1}^2 C_1^{jk}x_jx_k,$$

one observes  as noticed in \cite{GS3}, that, assuming that $(\psi,\Phi)$ or equivalently $(\psi,\phi)$ is a solution of \eqref {ZS} (resp. \eqref{ZSbis}), then $(u,\rho)$ defined by

$$u(x,t)=\frac{1}{t}e^{i\frac{q(x)}{4t}}\bar{\psi} \left(\frac{x}{t},\frac{1}{t}\right),$$

$$\rho_j(x,t)=\frac{1}{t}\phi_j\left(\frac{x}{t},\frac{1}{t}\right)$$

is also a solution. \footnote{This result can be extended to the three-dimensional Zakharov-Schulman systems, replacing \eqref{phitwo} by $ L_2\phi_j=\frac{\partial}{\partial x_j} |\psi|^{4/3}, j=1,2,3. $}

This implies that a finite time blow-up occurs provided there exists a {\it localized} solitary wave  solution of the form $e^{-i\omega t}\psi (x)$ (see Corollary 3.1 in \cite{GS3} for Zakharov-Schulman systems in dimension two and three). On top of the focusing integrable DS II system we analyzed above, this can be the case for elliptic-elliptic Davey-Stewartson systems. We refer to \cite{Ci, Ci2, Ot1,Ot2,Ot3} for theoretical issues on ground state solutions of elliptic-elliptic DS systems and the possible associated blow-up and   to \cite{PSSW} for numerical analysis of the blow-up. We again emphasize that a refined analysis {\it \`a la Merle-Rapha\"{e}l \cite{MeRa}} is still missing for elliptic-elliptic DS systems.

\subsection{Solitary waves}
It is well known that hyperbolic NLS equations such as \eqref{HNLS} do not 
possess solitary waves of the form $e^{i\omega t}\psi (x)$ where $\psi$ is localized (see \cite {GS1}).

It was proven in \cite{GS1} that  non trivial solitary waves may exist for DS II type systems only when $\rho =1$ (focusing case) and $\beta\in (0,2).$ Note that the (focusing) integrable case corresponds to $\beta =1.$ Moreover  solitary waves with  {\it radial} (up to translation) profiles can exist only when $\rho=1$ and $\beta=1$, that is in the focusing integrable case.

Those results  suggest that solitary waves for the focusing DS II systems exist only in the integrable case and this might be due to the new symmetry of the system we were alluding to above in this case.

To summarize, one is led to conjecture that neither the existence of 
the lump nor the associated Ozawa blow-up persist in the focusing DS II non integrable case.

\subsection{Line solitary waves}
For solutions depending only on $x,$ \eqref{DSII} with $\epsilon =1$ reduces to the one-dimensional NLS equation

\begin{equation}\label{1D}
i\partial_t\psi+\partial_{xx}\psi+2\rho(1-2\beta)|\psi|^2\psi=0.
\end{equation}

When $\alpha=2\rho(1-2\beta)>0,$ (in particular in the integrable 
focusing case $\beta=1, \rho=-1$ of DS II), \eqref {1D} is focusing and 
possesses the explicit solitary wave 
$\psi(x,t)=e^{it}\sqrt{\frac{2}{\alpha}} \frac{1}{\cosh x}$ 
which is asymptotically stable within \eqref{1D}. An interesting 
question is that of its {\it transversal} stability with respect of \eqref{DSII}. This issue has been investigated theoretically by Rousset and Tzvetkov for various nonlinear dispersive PDE's (\cite{RT1,RT2,RT3}), in particular for the cubic two-dimensional NLS, but not to our knowledge in the context of Davey-Stewartson systems. It is proven in \cite{PRKD} that the KdV line soliton is unstable for short  wave transverse perturbations in the context of the hyperbolic NLS. A numerical study for the hyperbolic NLS and the integrable DS II system can be found in \cite{MR}. We plan to come back to theoretical and numerical issues for the non integrable DS II systems in a subsequent work.

On the other hand, an explicit formula is given in \cite{FPS} for the 
interaction of an N-line soliton and a lump of the integrable focusing DS II system. Since the lump does not seem to persist in the non integrable case, this situation has probably no counterpart then.

\section{Numerical methods}\label{sec:num}

In this section, we discuss the numerical methods used to 
compute the time-evolution of the solution and in particular 
how to identify the type of blow-up in case certain norms of the 
solution diverge. 

\subsection{Numerical methods for the time-evolution}\label{sec:numtime}

For the numerical integration of \eqref{DSII}, we use
a Fourier spectral method in $x$. The reasons for this choice are the 
excellent approximation properties of spectral methods for smooth 
functions, and the  minimal introduction of numerical 
dissipation which (in principle) could overwhelm the dispersive 
effects of DS II
we want to study. 

The discretization in Fourier space leads to a system of 
(stiff) ordinary differential equations for the Fourier 
coefficients of $\psi$ of the form
\begin{equation}
   \partial_t \widehat{\psi}=\mathcal{L}\widehat{\psi}+\mathcal{N}(\psi),
    \label{fnlsfourier}
\end{equation}
where $\mathcal{L}=-i\epsilon(\xi_{1}^{2}-\xi_{2}^{2})$, and where
$\mathcal{N}(\psi)=2i\rho/\epsilon\mathcal{F}\left( 
\lbrack(\Delta^{-1}(\partial_{yy}+(1-2\beta\partial_{xx})) |\psi|^2\rbrack\psi\right)$ 
denotes the nonlinearity. It is an advantage of Fourier methods 
that the derivatives and thus the operator $\mathcal{L}$ are 
diagonal. For equations of 
the form (\ref{fnlsfourier}) with diagonal $\mathcal{L}$, there are many efficient high-order time 
integrators. For DS II
the performance of several fourth order methods was recently 
compared in \cite{KR2}. 
It was shown that in the focusing case a {\it composite Runge-Kutta method} 
\cite{D} performed best, which we will also use in the defocusing 
case. 

The numerical precision is 
controlled via the numerically computed energy (\ref{energy})
 Due to 
unavoidable numerical errors, the computed energy will depend on 
time. It was shown in \cite{KR} that the quantity 
$\Delta_{E}=E(t)/E(0)-1$ can be used as an indicator 
of the numerical accuracy that overestimates the $L_{\infty}$ norm of 
the difference between numerical and exact solution by roughly two 
orders of magnitude. We always aim at a $\Delta_{E}$ smaller than 
$10^{-4}$ which guarantees plotting accuracy, and typically we 
achieve $\Delta_{E}<10^{-6}$ or better. 

\subsection{Dynamical rescaling}
A useful tool in the numerical study of blow-up in NLS equations are dynamically 
rescaled codes, see for instance
\cite[Chapter 6]{SS99} and references therein. For the 
elliptic-elliptic DS, this was used in \cite{PSSW} to study numerically blow-up in this 
system and to show that the latter behaves essentially as the standard NLS. 
Using the same technique for the focusing DS II, we put 
\begin{equation}
    X = \frac{x}{L(t)},\quad Y = \frac{y}{L(t)},\quad 
    \tau=\int_{0}^{t}\frac{dt'}{L^{2}(t')},\quad \Psi(\xi,\eta,\tau) = 
    L(t)\psi(x,y,t)
        \label{resc}.
\end{equation}
Equation (\ref{DSIIbis}) implies for $\Psi$ 
\begin{equation}
    \label{DSIIresc}
\begin{array}{ccc}
i\epsilon 
\partial_{\tau}\Psi+\epsilon^{2}\partial_{XX}\Psi-\epsilon^{2}\partial_{YY}\Psi+i\epsilon a(X 
\partial_{X}\Psi+Y\partial_{Y}\Psi+\Psi)&&\nonumber\\
-2\left(\beta\Phi+\left|\Psi\right|^{2}\right)\Psi & = & 0,
\\
\partial_{XX}\Phi+\partial_{YY}\Phi+2\partial_{XX}\left|\Psi\right|^{2} & = & 0, 
\end{array}
\end{equation}
where $a = \partial_{\tau}\ln L$.
Under this rescaling, the $L^{2}$ norm of $\psi$ with respect to $x$ 
and $y$ is equal to the $L^{2}$ norm of $\Psi$ with respect to $X$ 
and $Y$. 

The scaling function $L$ can be chosen in a way that $L(\tau)$ 
vanishes for $\tau\to\infty$ corresponding to $t\to t^{*}$, the 
blow-up time. 
It is then expected at blow-up that both $a\to a^\infty$ and $\Psi 
\to \Psi^\infty$ become $\tau$-independent.  Equation 
(\ref{DSIIresc}) then becomes a PDE in $X$ and $Y$ only which would 
give the asymptotic profile of the selfsimilar blow-up. There is no 
reason to believe that this equation reduces to an ODE in generic 
situations, nor that it has radially symmetric solutions in this case. For $a^{\infty}=0$, a case which is not expected for 
blow-up in NLS solutions, these would be defining equations for 
solitary waves.

The choice of the scaling factor $L(t)$ is done for numerical 
convenience. Typically it is fixed by demanding that one of the norms 
diverging at blow-up for $\psi$ is kept constant, for instance the 
$L^{\infty}$ norm of $\Psi$. Numerically it is better though to keep an 
integral norm of $\Psi$ constant. Since the $L^{2}$ norm is anyway 
invariant, we choose the $L^{2}$ norm of the $X$-derivative. 
This leads to 
\begin{equation}\label{L}
L(t) = \frac{ | \partial_{X} \Psi(\tau,\cdot)|_{2}}{ | 
\partial_{x} \psi(t, \cdot)|_{2}}
\end{equation}
where  $| \partial_{X} \Psi|_{2}$ is chosen to be 
constant. We can read off the time-evolution of 
$L$ from (\ref{L}) and \eqref{resc} by differentiating the $L^{2}$ 
norm of $\partial_{X}\Psi$ with respect to $\tau$ and by using 
(\ref{DSIIresc}) to eliminate the $\tau$-derivatives, which leads after some partial integrations to
\begin{equation}
    a(\tau)=\frac{2}{\epsilon|\partial_{X} \Psi |_2^{2} 
    }\int_{\mathbb{R}^{2}}(\beta\Phi+ |\Psi|^{2})  
    \Im(\bar{\Psi}\partial_{XX}\Psi ) dX dY
    \label{L2x}.
\end{equation}

This allows us in principle to study the type of the 
blow-up for DS II in a similar way as it has been done for generalized Korteweg-de 
Vries equations in \cite{KP2013} by numerically integrating 
(\ref{DSIIresc}). 
But it was shown 
numerically in \cite{KP2013} that generic rapidly decreasing hump-like 
initial data lead to a 
tail of dispersive oscillations towards spatial infinity with slowly 
decreasing amplitude. Due to the imposed periodicity (in our numerical domain), these 
oscillations reappear after some time
on the opposing side of the computational domain 
and lead to numerical instabilities in the dynamically rescaled 
equation. The source of these problems are the terms 
$X\Psi_{X}+Y\Psi_{Y}$ in 
(\ref{DSIIresc}) since $X$, $Y$ are large at the boundaries of the 
computational domain, which has to be chosen large enough to allow 
for the `zooming in' effect due to a smaller and smaller $L$. Small numerical errors tend to be amplified by 
these terms.  For gKdV this could be addressed 
by using high resolution in time and large computational domains, a 
resolution which is difficult to achieve for $2+1$ dimensions. For 
DS II there is the additional problem of the \emph{modulational 
instability} of focusing NLS equations which shows in numerical 
computations in the form of an increase of the Fourier coefficients 
for the high wave numbers.  The consequence of 
this is that  we 
cannot compute long enough with the dynamically rescaled code 
to get conclusive results. Instead we 
integrate DS II directly, as described above, and then we use some 
post-processing to characterize the type of blow-up via the above 
rescaling (\ref{resc}) and (\ref{L2x}).

Under the hypothesis that $L(\tau) \sim\exp(-\kappa \tau)$ close to 
blow-up with $\kappa>0$ some positive constant, (\ref{resc}) yields a 
connection between $t$ and $\tau$,   
\begin{equation}
    L(t) \propto \sqrt{t^{*}-t}
    \label{Lt}.
\end{equation} 
With (\ref{L}) and (\ref{resc}), this 
implies
\begin{equation}
    | \partial_{x} \psi(t,\cdot) |_{2}^{2}\propto (t^{*}-t)^{-1},\quad
    |\psi(t,\cdot) |_{\infty}\propto 
    (t^{*}-t)^{-1/2}.    
    \label{genscal}
\end{equation}
In the mass critical case for  NLS, one finds a correction to \eqref{L} in the form 
\begin{equation}
L(t) \propto \sqrt{\frac{t^{*}-t}{\ln |\ln (t^{*}-t)|}},
    \label{L2scal}
\end{equation}
i.e., one has $\tau\propto \ln (t^{*}-t)(1-\ln |\ln (t^{*}-t)|)$ 
instead of $\tau\propto \ln (t^{*}-t)$. This so-called log-log-scaling regime for mass critical NLS has been rigorously proved in \cite{MR}.
It cannot be expected that  
logarithmic corrections can be seen in our simulations.

\subsection{Singularity tracing in the complex plane}

In the case of fractional NLS equations with cubic nonlinearity, it was observed in 
\cite{KSM} that the codes continue to run even if a finite-time 
blow-up is reached. This is in contrast to NLS equations with higher 
nonlinearity where the $L^{\infty}$ norm of 
the solution becomes so large that the computation of the nonlinear 
terms in the  equation leads to an overflow error. Since the 
nonlinearity of NLS is also cubic, we 
identify an appearing singularity as follows (see also \cite{KR2013a, KR2, SSF}):
Recall that in the complex plane, a (single) singularity $z_0\in \C$ of a real 
function $f$, such that $f(z) \sim (z-z_{0})^{\mu}$, with $\mu\not \in \mathbb Z$, results in the following asymptotic behavior for the corresponding Fourier transform
\begin{equation}
    |\widehat{f}(k)|\sim 
    \frac{1}{k^{\mu+1}} e^{-k\delta},\quad |k| \gg 1,
    \label{fourasymp}
\end{equation}
where $\delta=\text{Im}\, z_{0}$. The quantity $\mu$ thereby 
characterizes the type of the singularity.

In \cite{KR2013a, KR2}  this 
approach was used to quantitatively identify the time where the 
singularity hits the real axis, i.e., where the real solution becomes 
singular, since it was shown that the quantity $\delta$ can be reliably identified 
from a fitting of the Fourier coefficients. This is not true for 
$\mu$ though since the numerical inaccuracy is too large. 
In the case of focusing NLS, it was shown in \cite{KR2} that 
the best results are obtained when the code is stopped once the 
singularity is closer to the real axis than the minimal resolved 
distance via Fourier methods, i.e.,
\begin{equation}
    m:=2\pi \frac{D}{N},
    \label{mres}
\end{equation}
with $N\in \N$ being the number of 
 Fourier modes and $2\pi D$ the length of the computational domain in 
 physical space.  All values of $\delta <m$ 
 cannot be distinguished numerically from $0$. 
 
 Note that the 
 time at which the code is stopped because of the criterion above is \emph{not}
 the blow-up time itself. Rather, it is only the time where the code stops to be 
 reliable. The blow-up time will be determined from the numerical 
 data by fitting to the scalings given in the previous subsection.

\section{Numerical results}

We aim here to investigate numerically if the properties displayed by the integrable DS II system persist in the non-integrable case. Are there "generic"? In particular does blow-up occur in the non-integrable DS II systems and are the solutions of the non-integrable defocusing DS II (and of the hyperbolic NLS) globally defined and disperse as $t\to \infty?$

The cubic NLS equation in $2+1$ dimensions can have blow-up, as the elliptic-elliptic DS systems (\cite{GS, PSSW}). As previously recalled, 
DS II solutions can also have blow-up in the integrable case. Results by Sung 
\cite{Su4} establish for the integrable case global existence in time for 
initial data $\psi_{0}\in L_{p}$, $1\leq p < 2$ with a Fourier 
transform $\mathcal{F}[\psi_{0}]\in L_{1}\cap L_{\infty}$ subject to the 
smallness condition 
\begin{equation}
    |\mathcal{F}[\psi_{0}]|_1|\mathcal{F}[\psi_{0}]|_\infty
    <\frac{\pi^{3}}{2}\left(\frac{\sqrt{5}-1}{2}\right)^{2}
    \label{sungcond}
\end{equation}
in the 
focusing case. There is no such condition in the defocusing case. 
Corresponding results for the non-integrable cases are not yet known. 
Note that condition (\ref{sungcond}) has been established for the DS II 
equation with $\epsilon=1$. After the coordinate change $x'=x/\epsilon$, $y'=y/\epsilon$, 
$t'=t/\epsilon$  condition 
(\ref{sungcond}) takes for the initial data $\psi_{0}(x,y)=\exp(-x^{2}-y^{2})$
the form 
$$\frac{1}{\epsilon^{2}}\leq 
\frac{1}{8}\left(\frac{\sqrt{5}-1}{2}\right)^{2} \sim 0.0477.$$
This condition is not satisfied for the values of $\epsilon$ we study here. 
In \cite{KR} initial data of the form $\psi_{0}(x,y)=\exp(-x^{2}-\eta 
y^{2})$ were studied for the integrable case. No blow-up was observed 
for $\eta\neq 1$. In \cite{KR2} this was investigated in more detail, 
and it was found that there will be blow-up for the symmetric case 
$\eta=1$, but only in this case which will be studied in more detail. 


In the following we will always consider the initial data 
$\psi_{0}=\exp(-x^{2}-y^{2})$ for $\epsilon=0.1$. The computation is 
done for $(x,y)\in [-2\pi,2\pi]\times[-2\pi,2\pi]$ with 
$N_{x}=2^{12}$ or $N_{x}=2^{13}$ and $N_{y}=2^{12}$ Fourier modes. 
The Fourier coefficients decrease in this case to machine precision 
($10^{-16}$ here, which means that due to rounding errors 
values of $10^{-14}$ can be reached)
for the initial data.
In the cases where 
there is no blow-up, the maximum of the solution appears for $t<0.3$. 
We compute until $t=0.6$. 
We work with $N_{t}=4000$ time steps.

\subsection{Hyperbolic nonlinear Schr\"odinger equation}

\begin{figure}[htb!]
  \includegraphics[width=\textwidth]{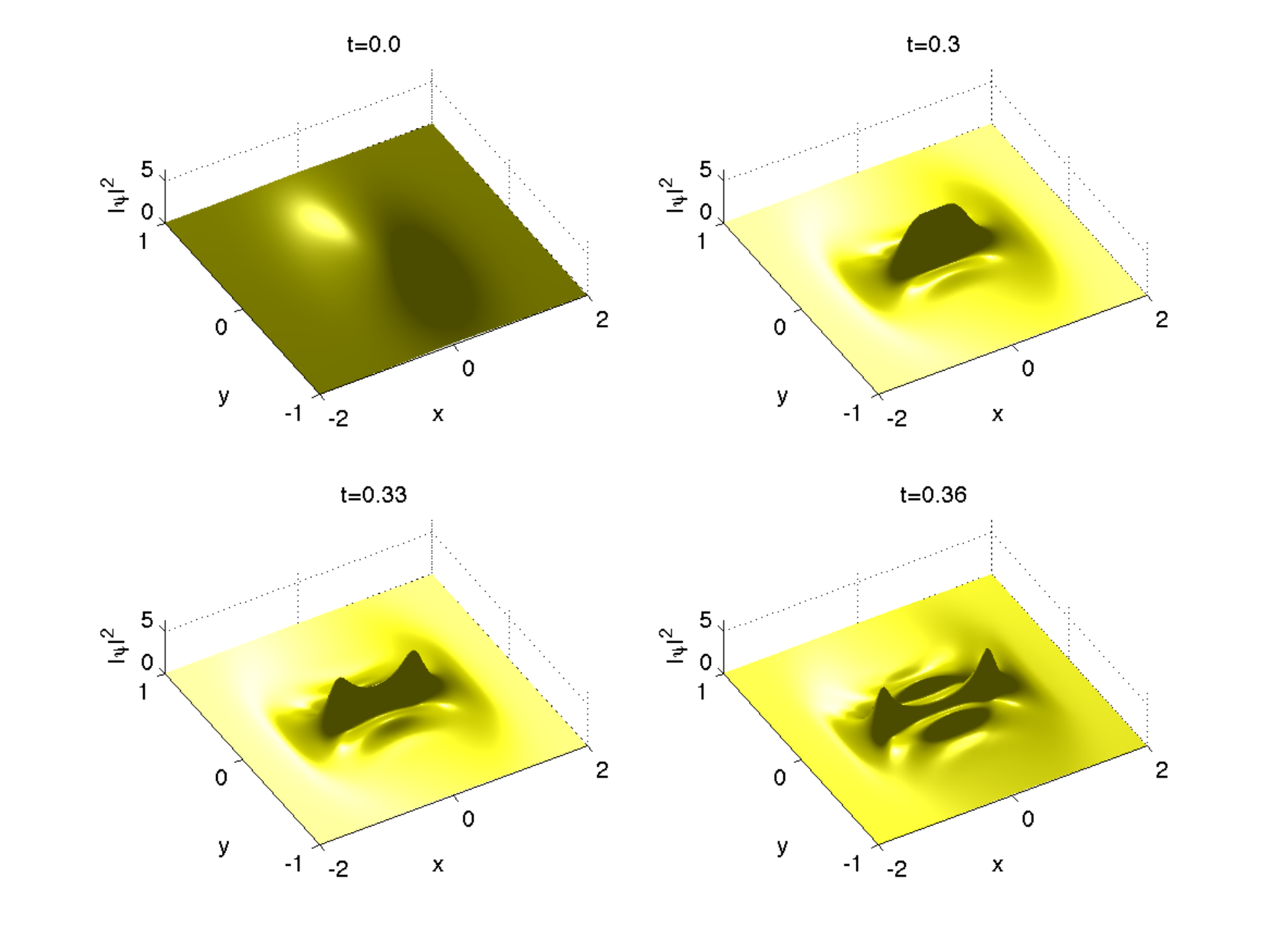}
 \caption{Solution to the hyperbolic NLS equation, equation 
 (\ref{DSII}) for $\beta=0$, for the initial 
 data $\psi_{0}=\exp(-x^{2}-y^{2})$ for $\epsilon=0.1$ at different 
 times.}
 \label{nlshyper4t}
\end{figure}

\begin{figure}[htb!]
  \includegraphics[width=0.7\textwidth]{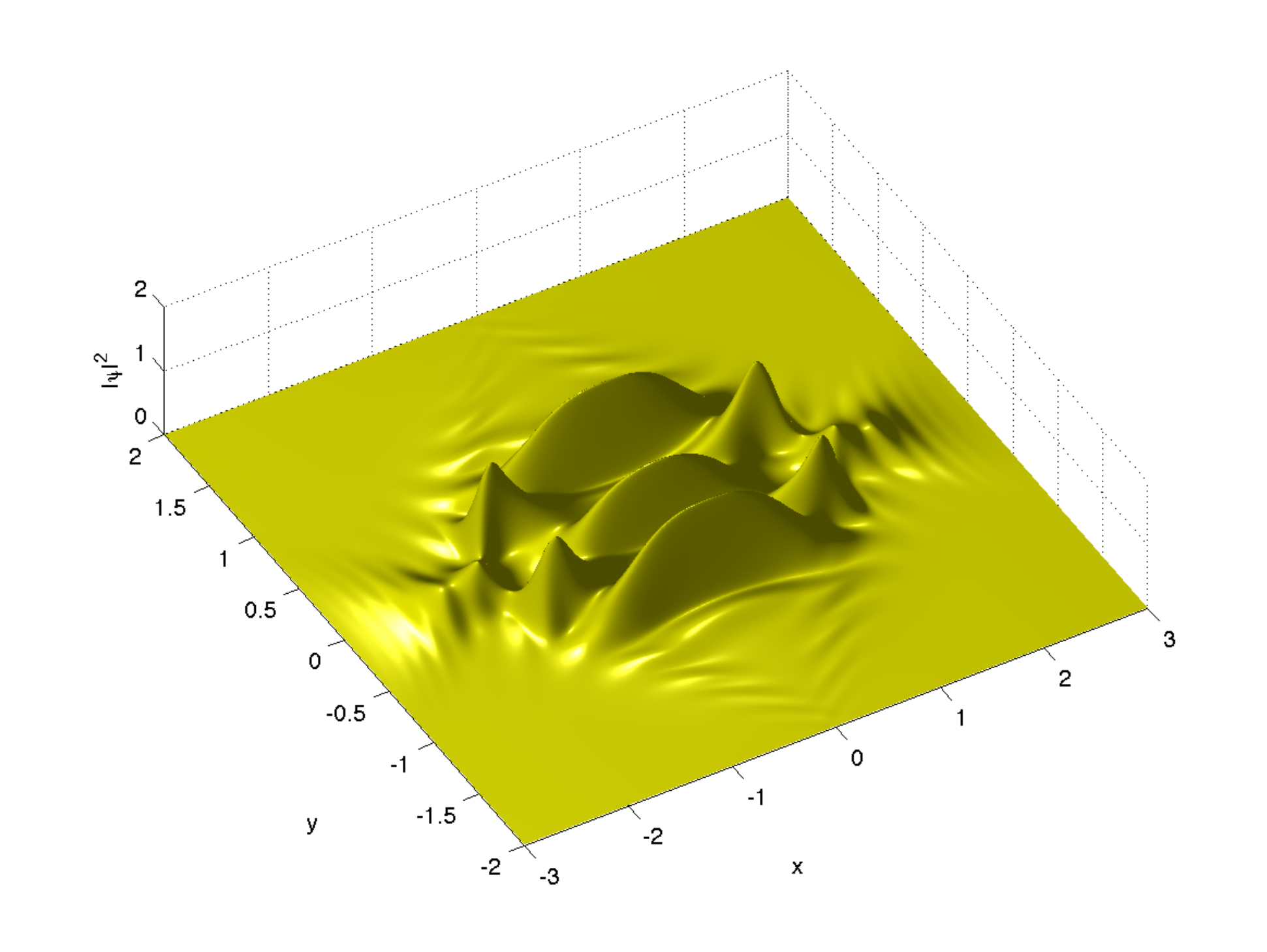}
 \caption{Solution to the hyperbolic NLS equation, equation 
 (\ref{DSII}) for $\beta=0$, for the initial 
 data $\psi_{0}=\exp(-x^{2}-y^{2})$ for $\epsilon=0.1$ at $t=0.6$.}
 \label{nlshyper06}
\end{figure}

We will first study the case of the hyperbolic NLS equation, i.e., 
$\beta=0$ in (\ref{DSII}). In Fig.~\ref{nlshyper4t} the solution for 
the initial data $\psi_{0}(x,y)=\exp(-x^{2}-y^{2})$ can be seen for 
different times. The initial pulse is visibly compressed in
the $y$-direction and defocused in the $x$-direction. At a given time, 
the initial hump decomposes into several smaller maxima.
At larger times, a regular pattern of peaks forms, see 
Fig.~\ref{nlshyper06}.

\begin{figure}[htb!]
  \includegraphics[width=0.45\textwidth]{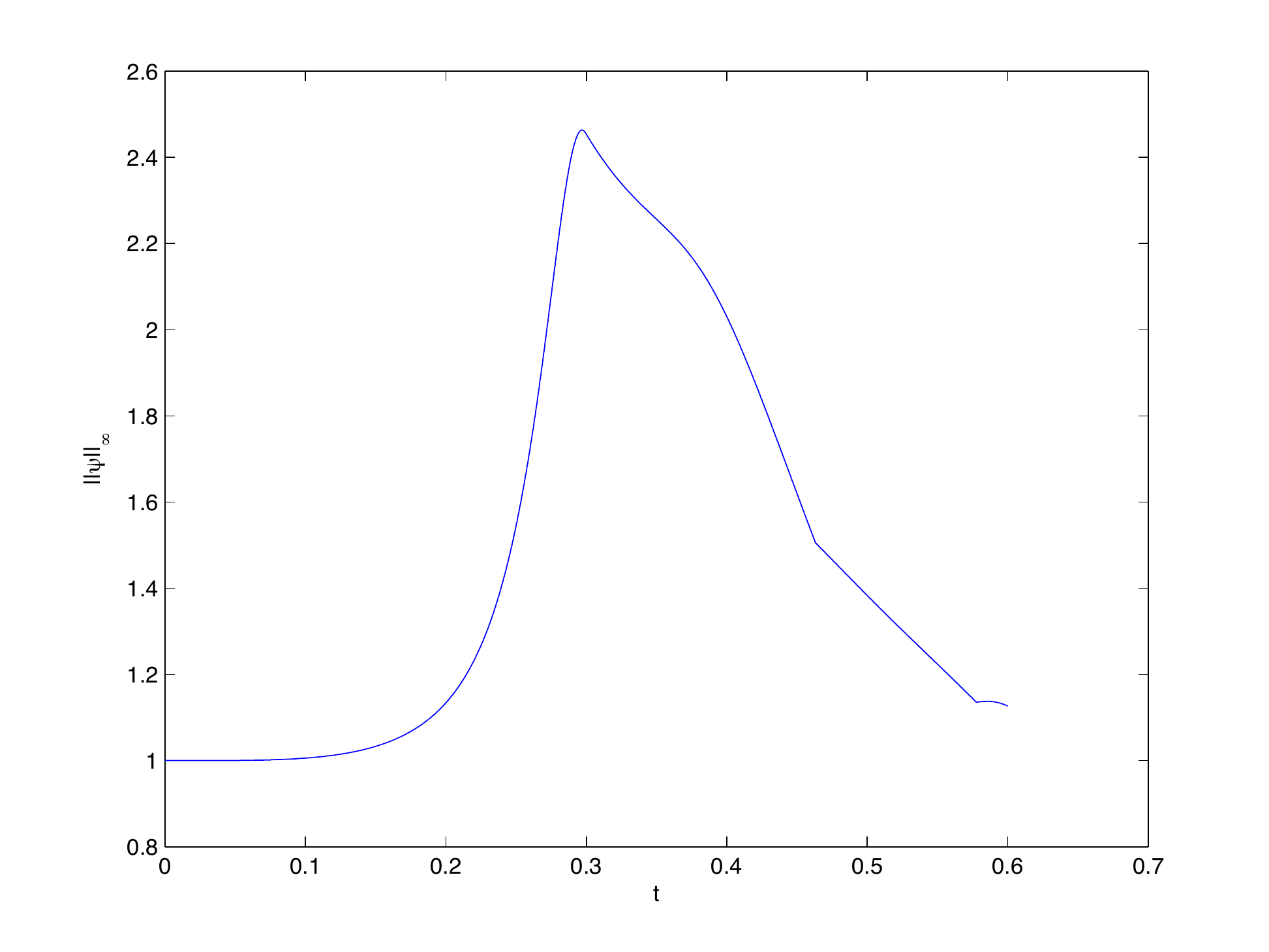}
  \includegraphics[width=0.45\textwidth]{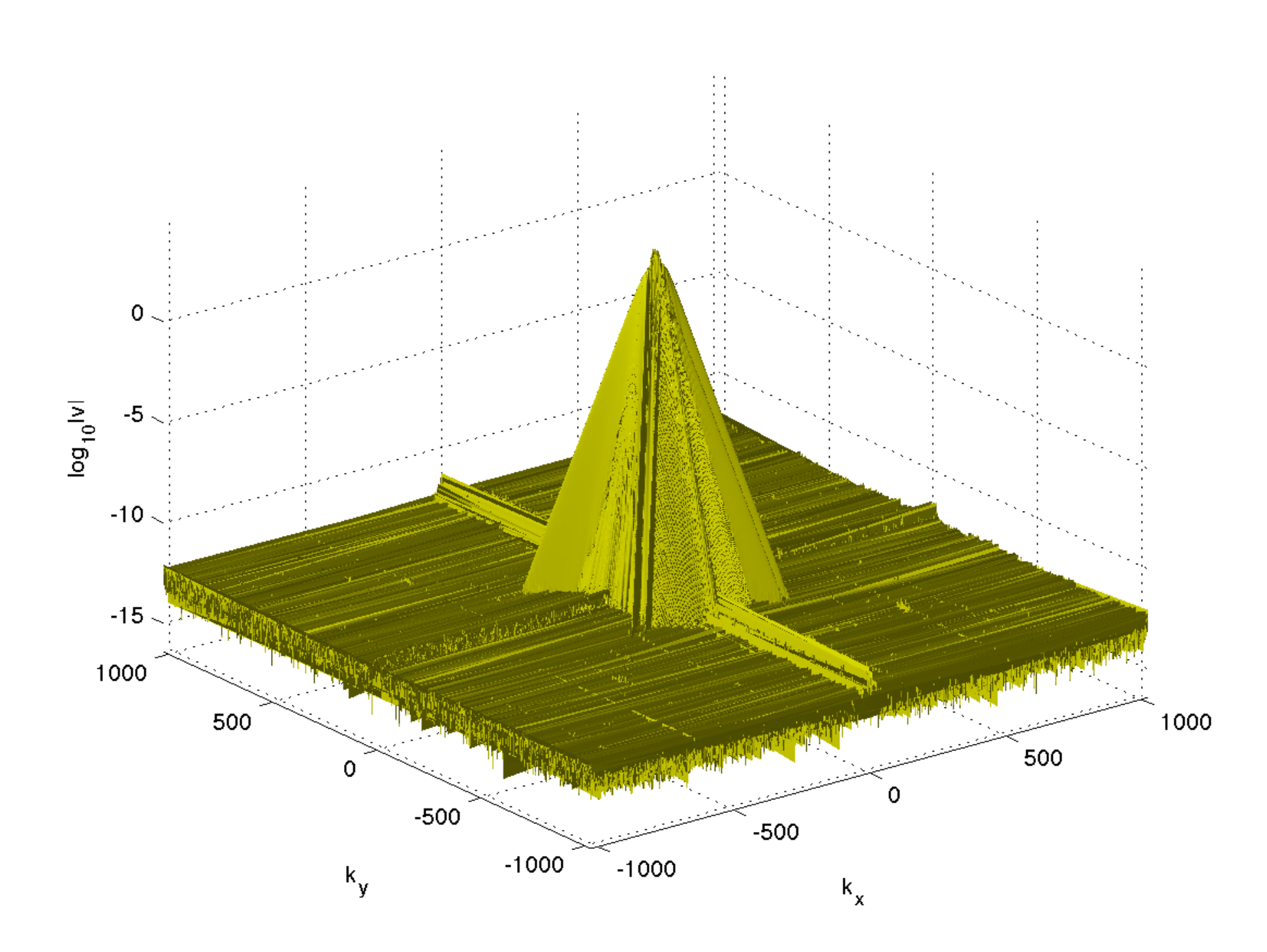}
 \caption{$L_{\infty}$-norm  on the left 
 and Fourier coefficients of the solution shown in 
 Fig.~\ref{nlshyper06} on the right.}
 \label{nlshyper06fourier}
\end{figure}

There is no indication of blow-up in this case as can be inferred 
from the $L_{\infty}$ norm of the solution in dependence of time in 
Fig.~\ref{nlshyper06fourier}. The norm increases until a time of 
roughly $0.3$ and decreases then monotonically. 
The solution is well resolved in Fourier space as is obvious from 
the same figure where the Fourier coefficients are 
shown. 

\subsection{Davey-Stewartson II solutions in the  focusing case 
$\rho=-1$}
\begin{figure}[htb!]
  \includegraphics[width=\textwidth]{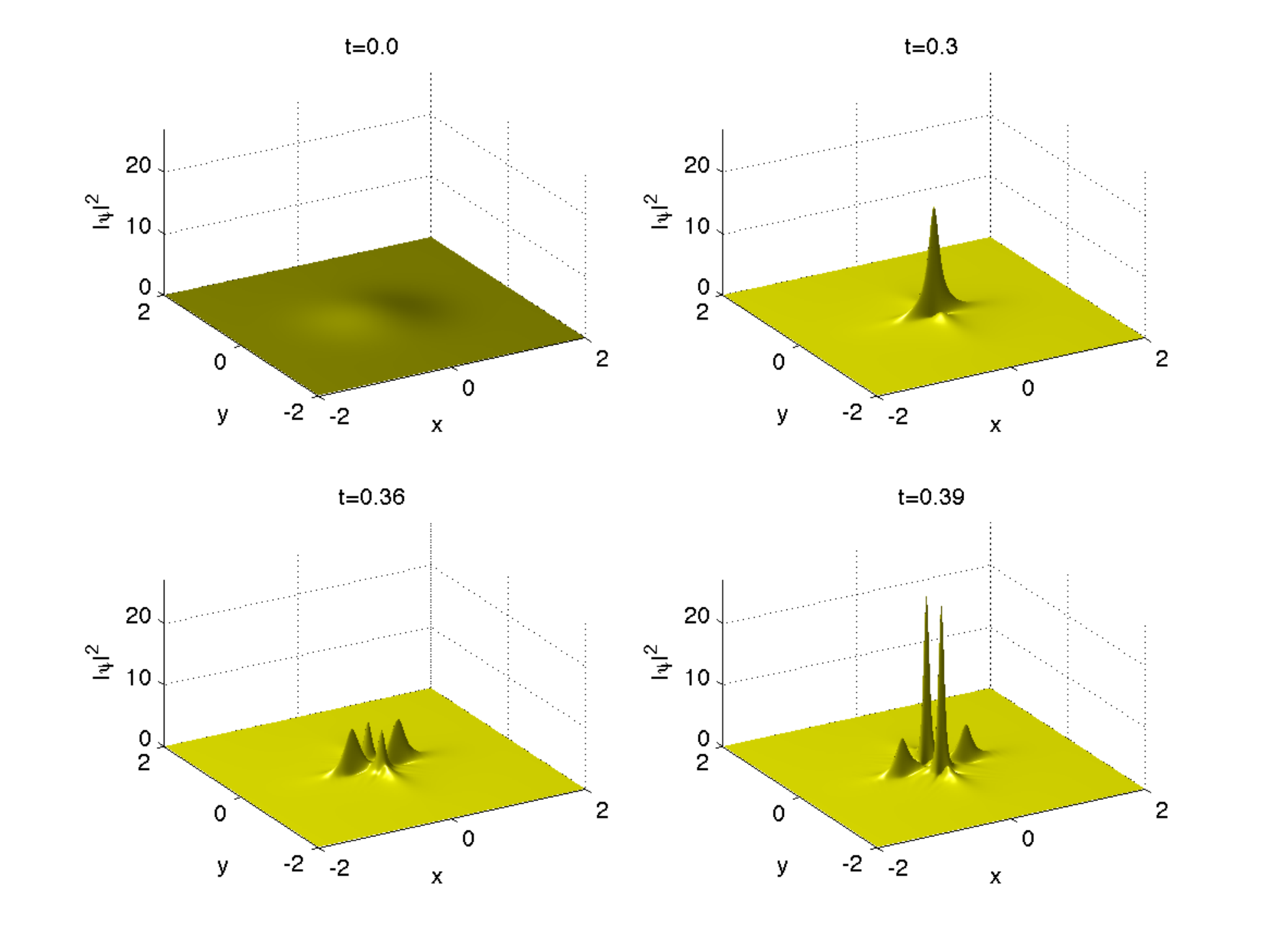}
 \caption{Solution to the DS II equation
 (\ref{DSII}) for $\beta=0.9$, for the initial 
 data $\psi_{0}=\exp(-x^{2}-y^{2})$ for $\epsilon=0.1$ at different 
 times.}
 \label{DSbeta09e014t}
\end{figure}

For $\beta<1$ the behavior of solutions to the focusing DS II system 
is similar to the hyperbolic NLS case, but the solutions get more and 
more focused. This can be seen for instance in 
Fig.~\ref{DSbeta09e014t} where the solution for the same initial 
data as before is shown for $\beta=0.9$ for several values of $t$. 
This indicates already that one effect of the nonlocal term $\Phi$ is to 
increase the focusing effect for $\rho=-1$ in the spatial 
direction which is for the hyperbolic NLS defocusing. It is clear 
that the peaks are much more pronounced in this case, but the main 
focusing effect is still in the $y$-direction.

\begin{figure}[htb!]
  \includegraphics[width=0.7\textwidth]{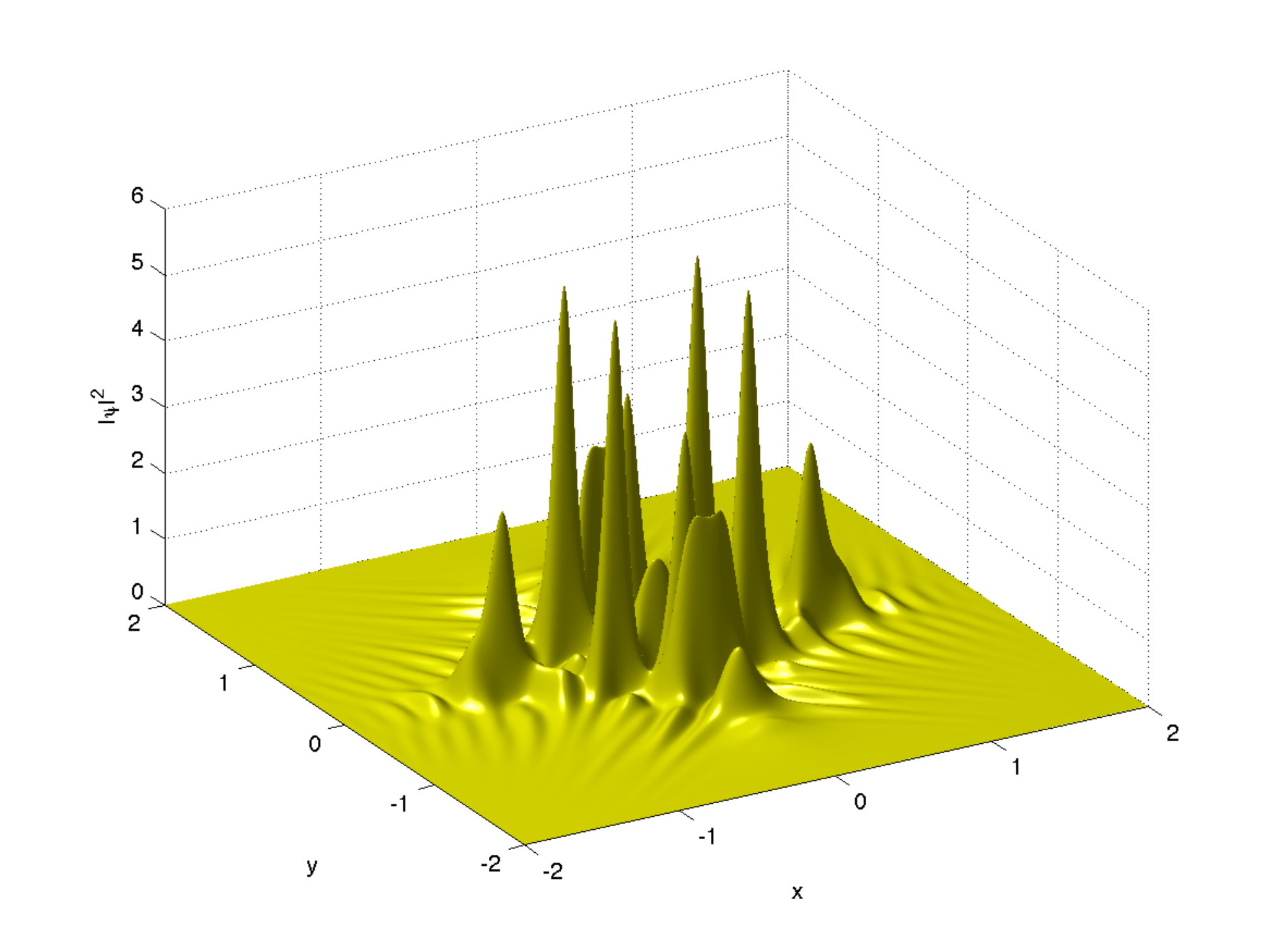}
 \caption{Solution to the DS II equation
 (\ref{DSII}) for $\beta=0.9$, for the initial 
 data $\psi_{0}=\exp(-x^{2}-y^{2})$ for $\epsilon=0.1$ at $t=0.6$.}
 \label{DSbeta09e01t06}
\end{figure}

At time $t=0.6$, one finds again a regular pattern of peaks and no 
indication of blow-up, see Fig.~\ref{DSbeta09e01t06}.

\begin{figure}[htb!]
  \includegraphics[width=0.45\textwidth]{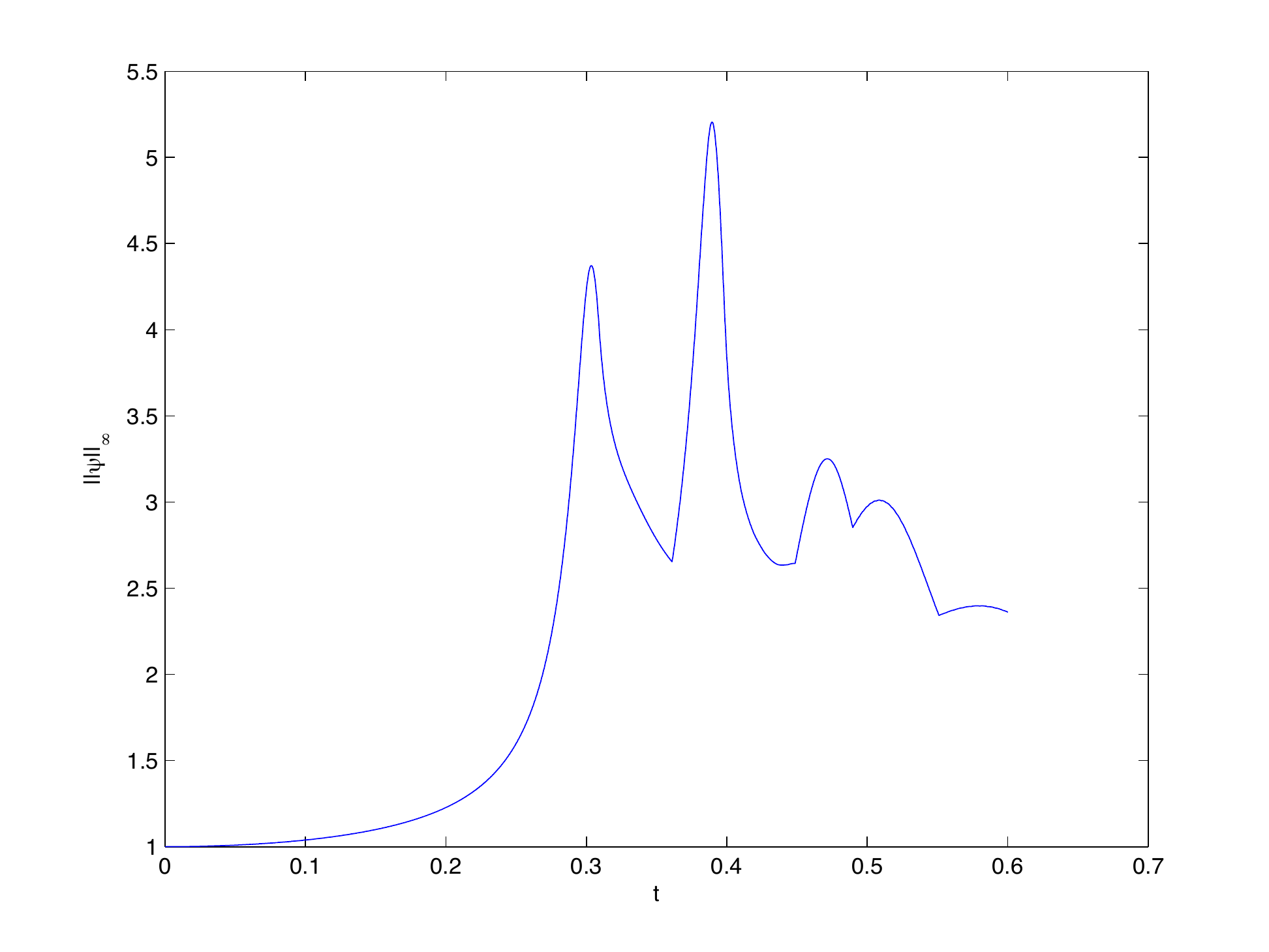}
  \includegraphics[width=0.45\textwidth]{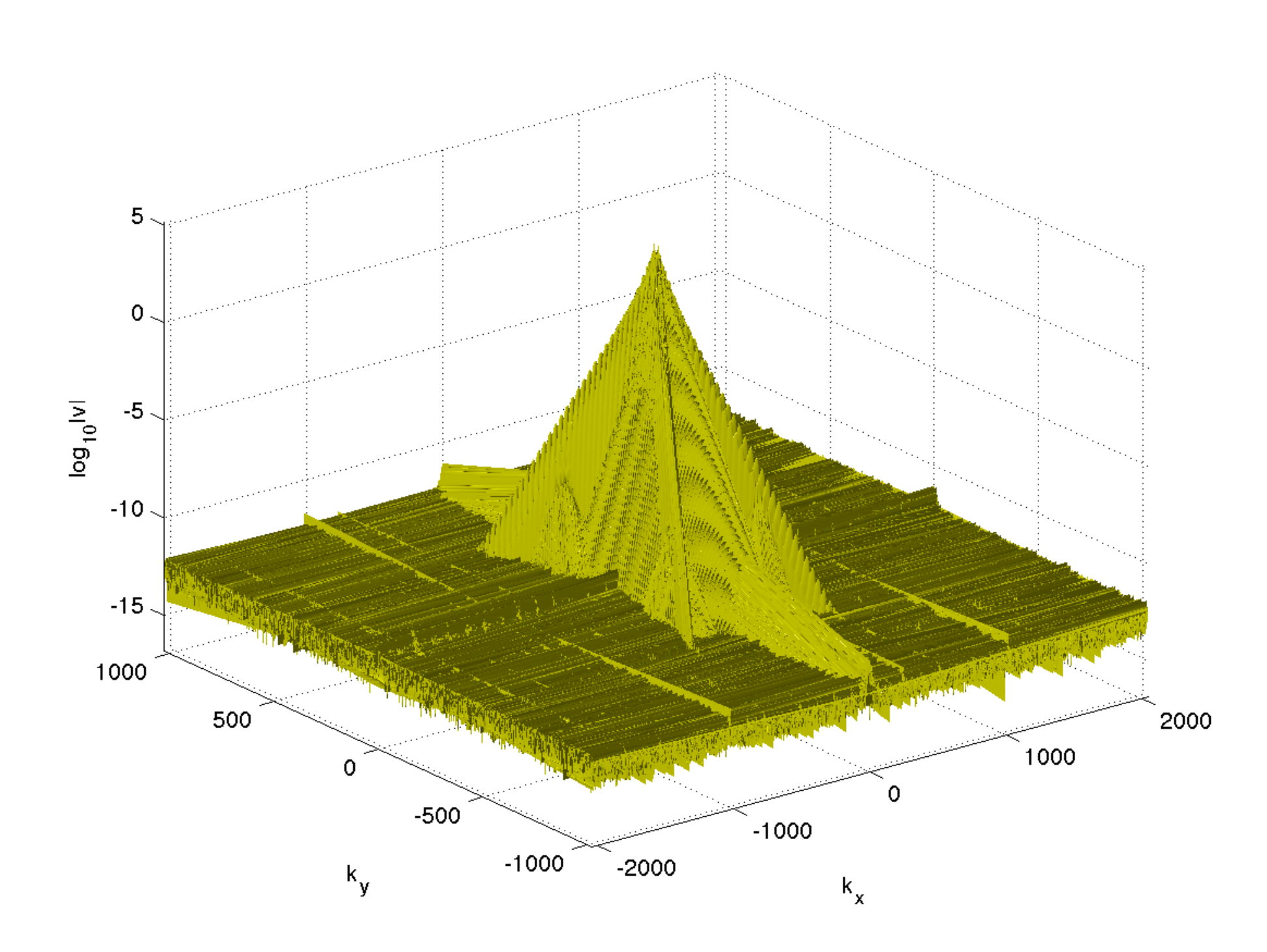}
 \caption{$L_{\infty}$-norm  on the left 
 and Fourier coefficients of the solution shown in 
 Fig.~\ref{DSbeta09e01t06} on the right.}
 \label{DSbeta09e01t06fourier}
\end{figure}

The appearance of several peaks essentially rules out the possibility 
of blow-up since this would be expected for symmetric initial data at 
$x=y=0$. If the central peak splits into several smaller peaks, it is 
not clear how this could lead to a blow-up at a later time. 
The $L_{\infty}$ norm of the solution in dependence of time in 
Fig.~\ref{DSbeta09e01t06fourier} shows that there are several peaks 
appearing also here in contrast to the hyperbolic NLS case, but that 
it appears to decrease with time after the highest maximum. 
The solution is again  well resolved in Fourier space as is visible in  
the same figure, but we had to use higher resolution and just reach machine 
precision in the $k_{y}$-direction. 

\begin{figure}[htb!]
  \includegraphics[width=0.7\textwidth]{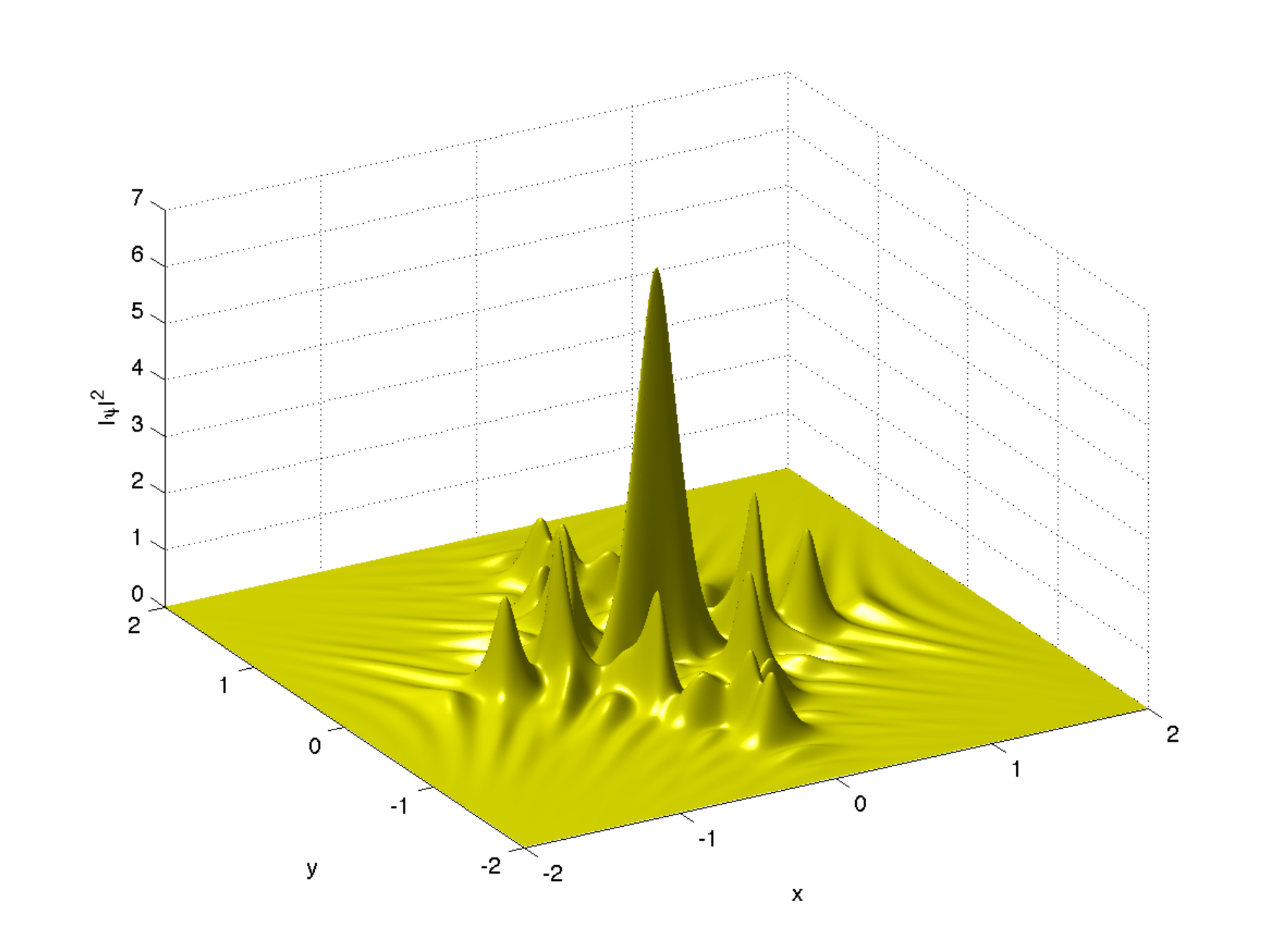}
 \caption{Solution to the DS II equation
 (\ref{DSII}) for $\beta=1.1$, for the initial 
 data $\psi_{0}=\exp(-x^{2}-y^{2})$ for $\epsilon=0.1$ at $t=0.6$.}
 \label{DSbeta11e01t06}
\end{figure}

It is remarkable that for the studied 
initial data, we obtain blow-up only for $\beta=1$ which will be 
discussed in detail in the next subsection. For $\beta>1$ we 
again observe a \emph{dispersive shock} as for $\beta<1$. This can be seen 
for $\beta=1.1$ in Fig.~\ref{DSbeta11e01t06} where no indication of 
blow-up can be seen.

\begin{figure}[htb!]
  \includegraphics[width=0.7\textwidth]{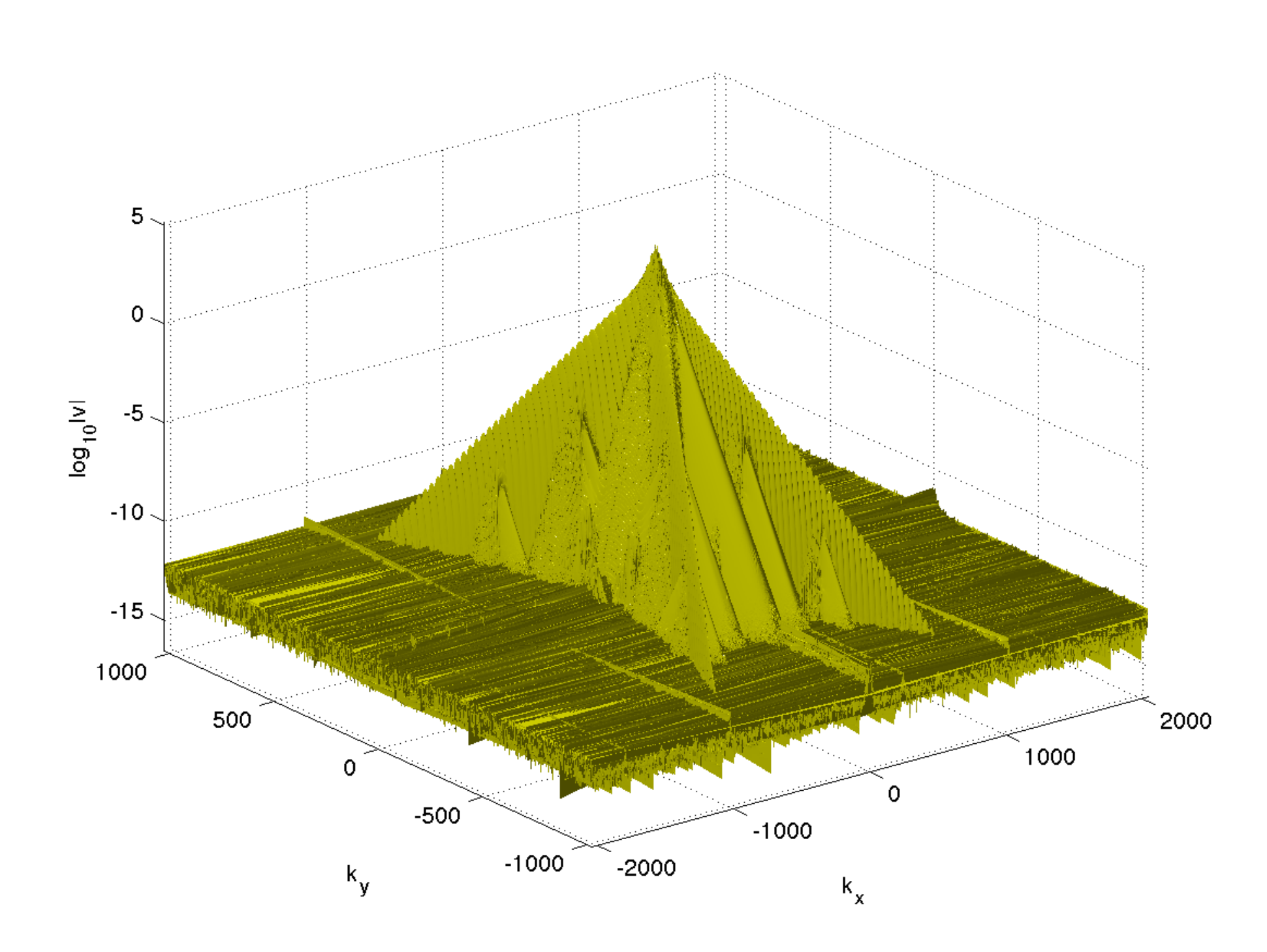}
 \caption{Fourier coefficients of the solution shown in 
 Fig.~\ref{DSbeta11e01t06}.}
 \label{DSbeta11e01t06fourier}
\end{figure}

The solution is also well resolved in Fourier space as is shown in 
Fig.~\ref{DSbeta11e01t06fourier}.

\begin{figure}[htb!]
  \includegraphics[width=0.7\textwidth]{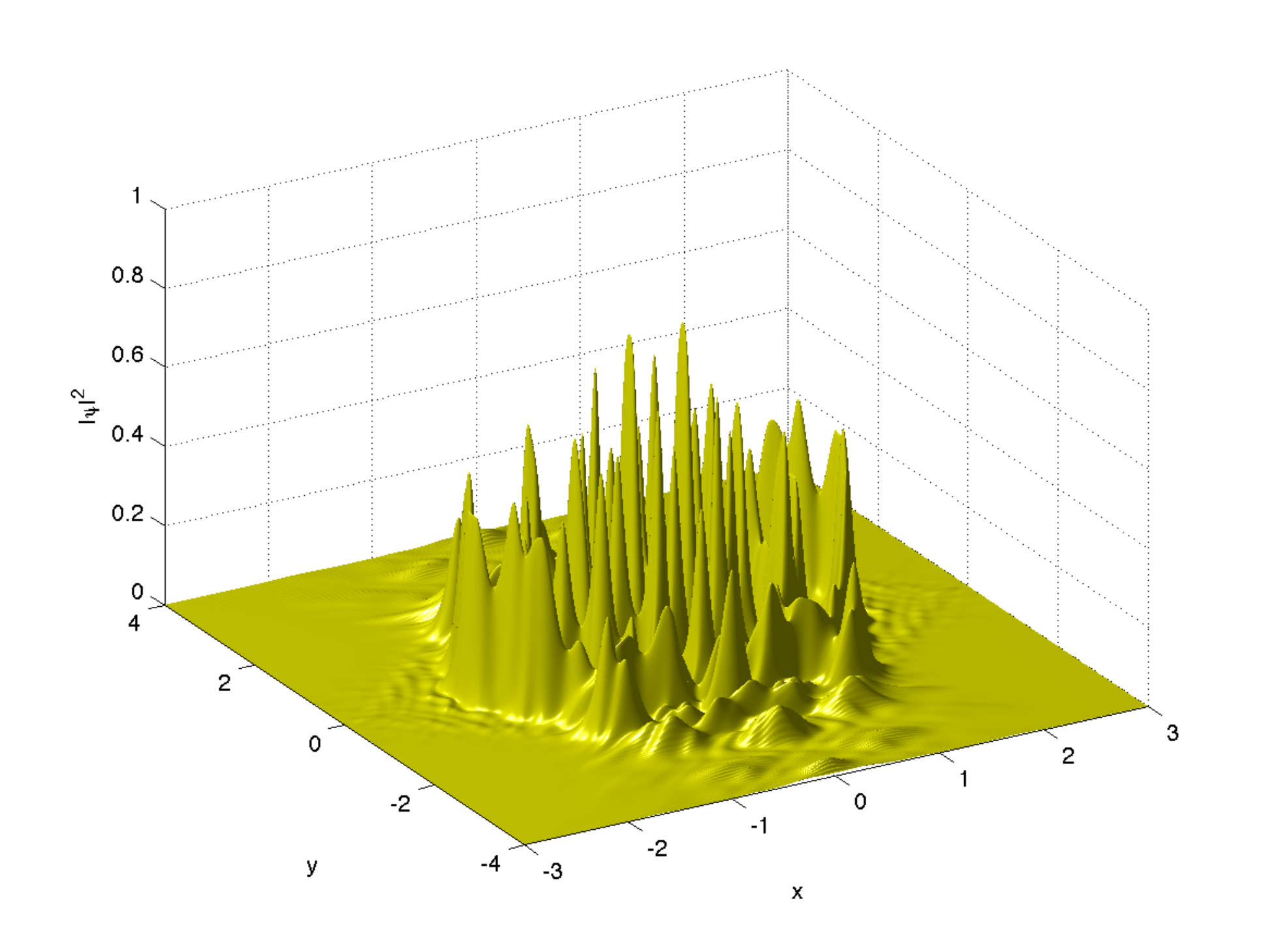}
 \caption{Solution to the DS II equation
 (\ref{DSII}) for $\beta=2$, for the initial 
 data $\psi_{0}=\exp(-x^{2}-y^{2})$ for $\epsilon=0.1$ at $t=0.6$.}
 \label{DSbeta20e01t06}
\end{figure}

For larger values of $\beta$ the focusing effect in $x$-direction 
becomes even stronger. The resulting dispersive shock in 
Fig.~\ref{DSbeta20e01t06} shows much more peaks than for smaller 
$\beta$, but there is again no indication of blow-up.

\begin{figure}[htb!]
  \includegraphics[width=0.7\textwidth]{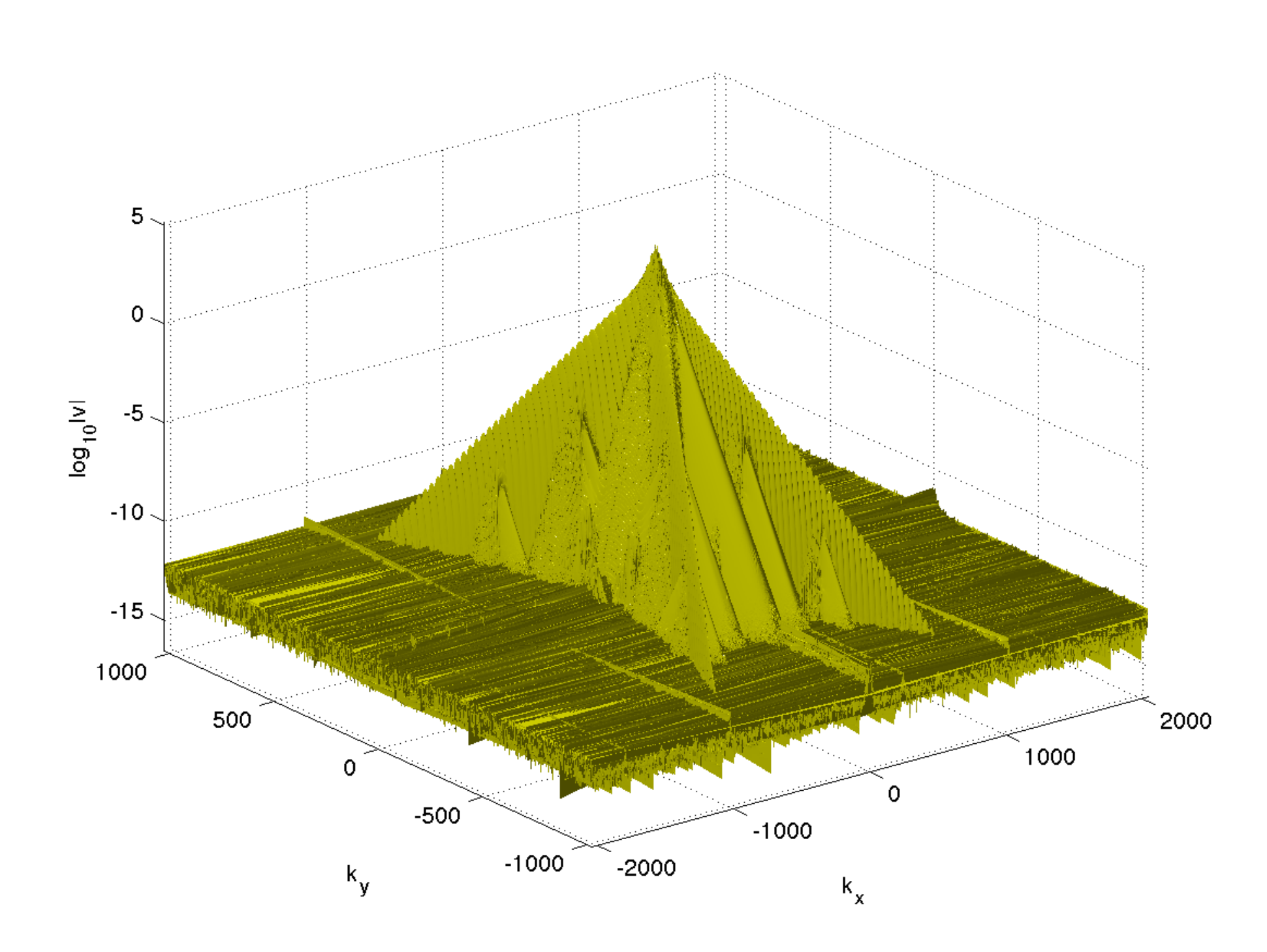}
 \caption{Fourier coefficients of the solution shown in 
 Fig.~\ref{DSbeta20e01t06}.}
 \label{DSbeta20e01t06fourier}
\end{figure}

\begin{figure}[htb!]
  \includegraphics[width=0.45\textwidth]{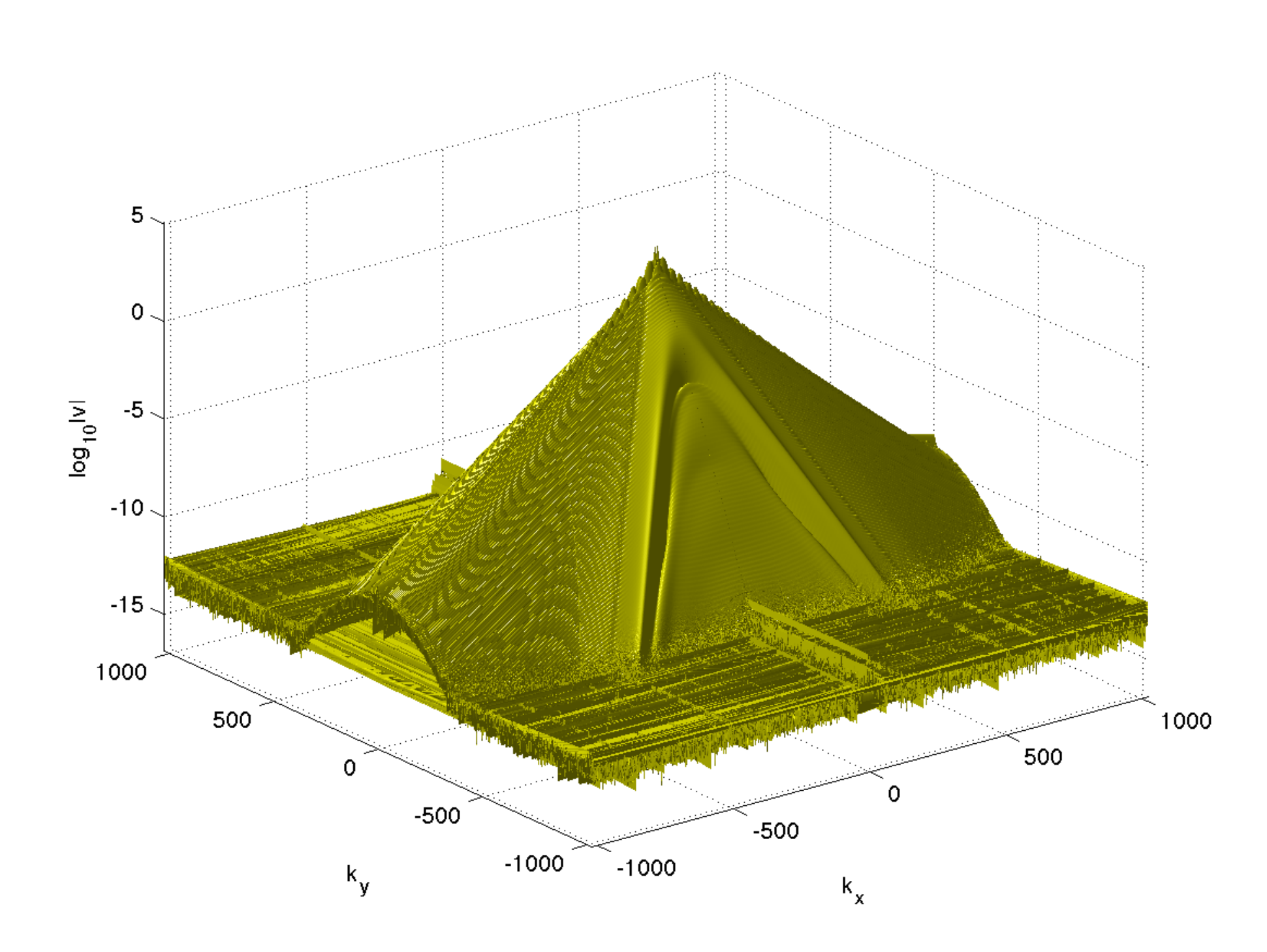}
  \includegraphics[width=0.45\textwidth]{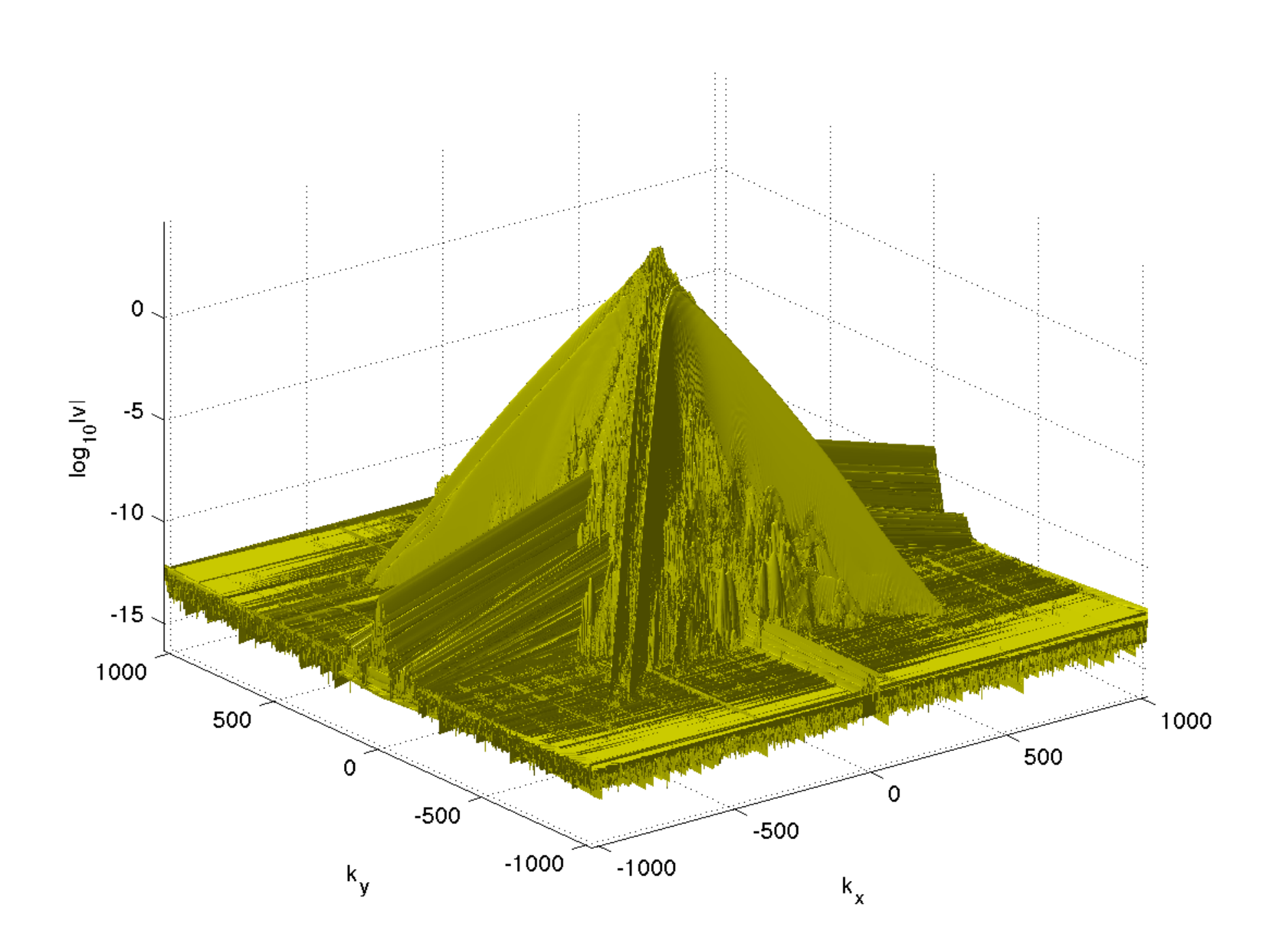}
 \caption{Fourier coefficients of the solution to the DS II equation 
 (\ref{DSII}) with $\beta=2$ for the initial 
 data $\psi_{0}=\exp(-x^{2}-y^{2})$ for $\epsilon=0.1$ at $t=0.18$ 
 on the left and at $t=0.6$ on the right for $N_{x}=N_{y}=2^{12}$.}
 \label{DSbeta20e01t06fourier2}
\end{figure}

In 
Fig.~\ref{DSbeta20e01t06fourier} we show that the situation of 
Fig.~\ref{DSbeta20e01t06} is numerically well resolved.
One problem in the numerical study of the semiclassical limit and 
possible blow-up of the focusing DS II equation is the well known 
modulational instability of the focusing NLS equation. This implies 
in the present context that a lack of spatial resolution leads to a spurious 
growing in time of the Fourier coefficients for the high wave numbers 
as was for instance discussed in \cite{etna}. Thus it is important to 
resolve well the maximum of the solution even if one only wants to 
study the situation at a later time for which less resolution is 
required. Such a situation can be seen in 
Fig.~\ref{DSbeta20e01t06fourier2}. The largest peak is observed at 
$t\sim 0.15$. At this time the peak is not resolved in 
$k_{x}$-direction up to machine precision. This leads to some 
artifacts in the Fourier coefficients at later times for the high 
wave numbers. Note, however, that the coefficients still decrease to 
$10^{-10}$ and that the solution is thus well resolved.  

\subsection{Blow-up in focusing DS II solutions}
\begin{figure}[htb!]
  \includegraphics[width=0.49\textwidth]{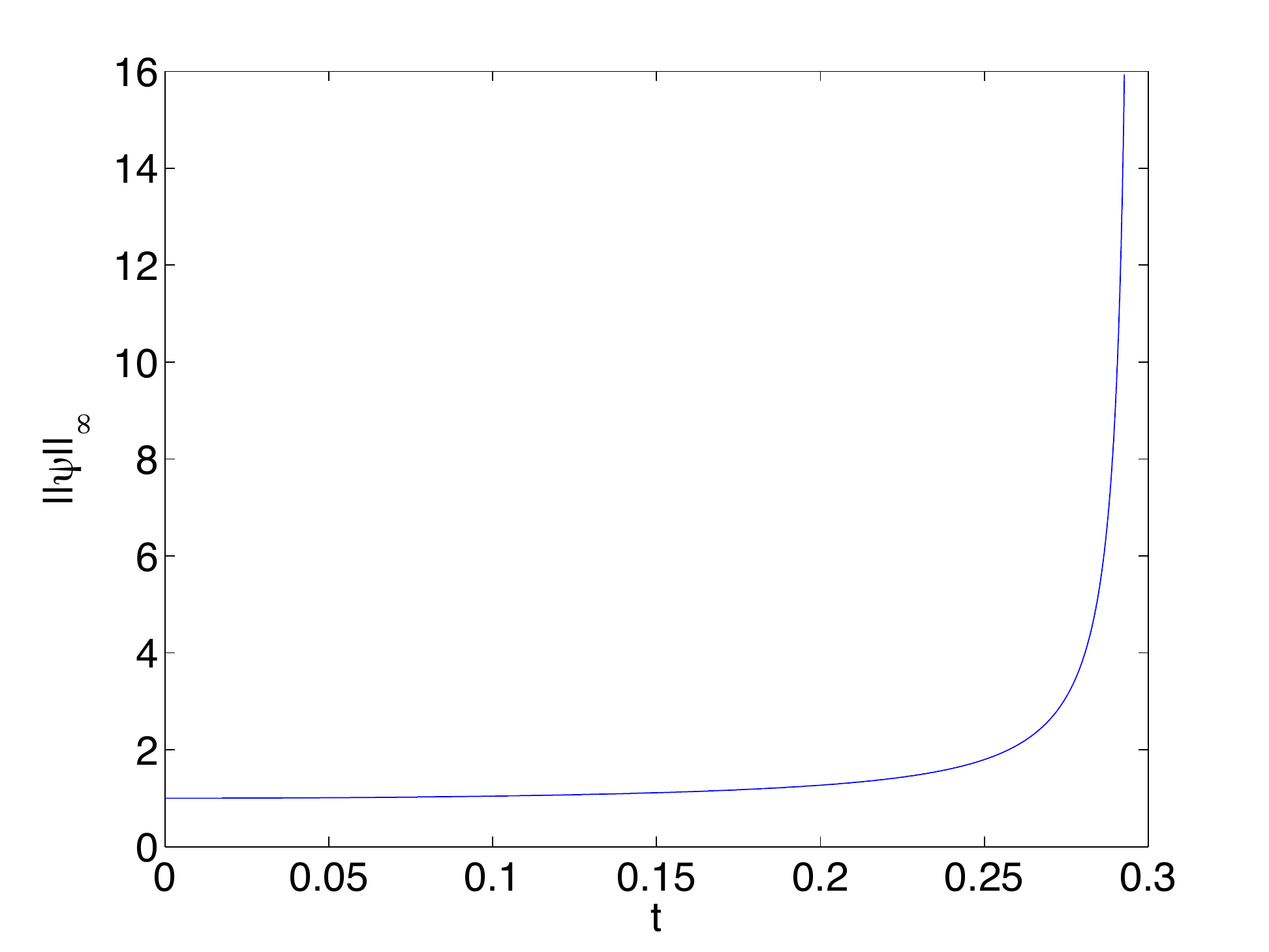}
  \includegraphics[width=0.49\textwidth]{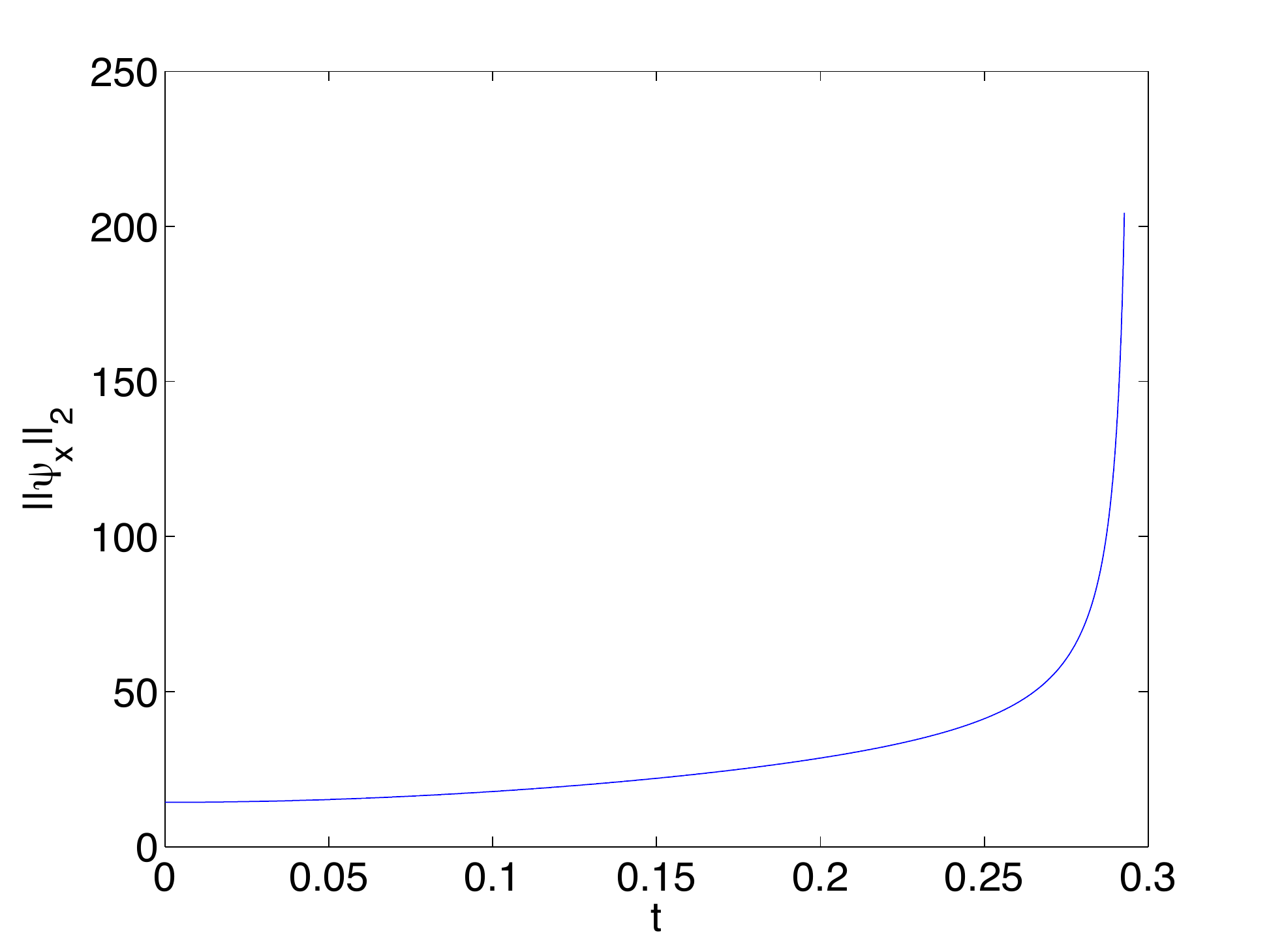}
 \caption{$L_{\infty}$-norm  of the solution to the DS II equation
 (\ref{DSII}) for $\beta=1$, for the initial 
 data $\psi_{0}=\exp(-x^{2}-y^{2})$ for $\epsilon=0.1$  in dependence 
 of $t$ on the left, and the corresponding $L_{2}$ norm of the 
 $\psi_{x}$ on the right.}
 \label{DSbeta1e01t06max}
\end{figure}

As was shown in \cite{KR2},  for initial data of the form 
$\psi_{0}=\exp(-(x^{2}+\gamma y^{2}))$ with $\gamma>0$, a dispersive 
shock is observed as in the previous subsection if $\gamma\neq1$.   
The situation is completely different in 
the integrable case $\beta=1$ for initial data 
$\psi_{0}=\exp(-x^{2}-y^{2})$ symmetric with respect to $x$ and $y$. 
In this case both the initial data and the equation are symmetric (up 
to complex conjugation) with respect to and exchange of the spatial 
coordinates. 
We use $N_{x}=N_{y}=2^{12}$ Fourier modes and $N_{t}=10^{4}$ time 
steps for $t<0.3$.
The initial hump is focused in both directions $x$ 
and $y$ in a symmetric manner. It turns out that the initial maximum 
grows without bounds. The code is stopped at $t=0.2927$ since the 
fitting of the Fourier coefficients  to the asymptotic formula 
(\ref{fourasymp}) indicates that a singularity is closer to the real 
axis than the minimal resolved spatial distance (\ref{mres}) (the 
fitting is done just in $x$-direction because of the symmetry with 
respect to $x\to y$). Note that the factor $\mu$ in 
(\ref{fourasymp}), which is in \cite{KR2} slightly smaller than 1, is 
here still larger than 1 which would indicate a cusp with finite 
values of the $L_{\infty}$ norm. This simply indicates that we do not 
have enough resolution for the blow-up. To obtain better 
approximations for the value of $\mu$, obviously higher resolution  on parallel 
computers would be needed as in \cite{KR2}, but we can make reliable
statements on the type of blow-up below 
with the serial computing used in the present paper. 
The 
$L_{\infty}$ norm of the solution and the $L_{2}$ norm of the 
$x$-derivative can be seen in 
Fig.~\ref{DSbeta1e01t06max}. Both seem to indicate a blow-up.

\begin{figure}[htb!]
   \includegraphics[width=0.49\textwidth]{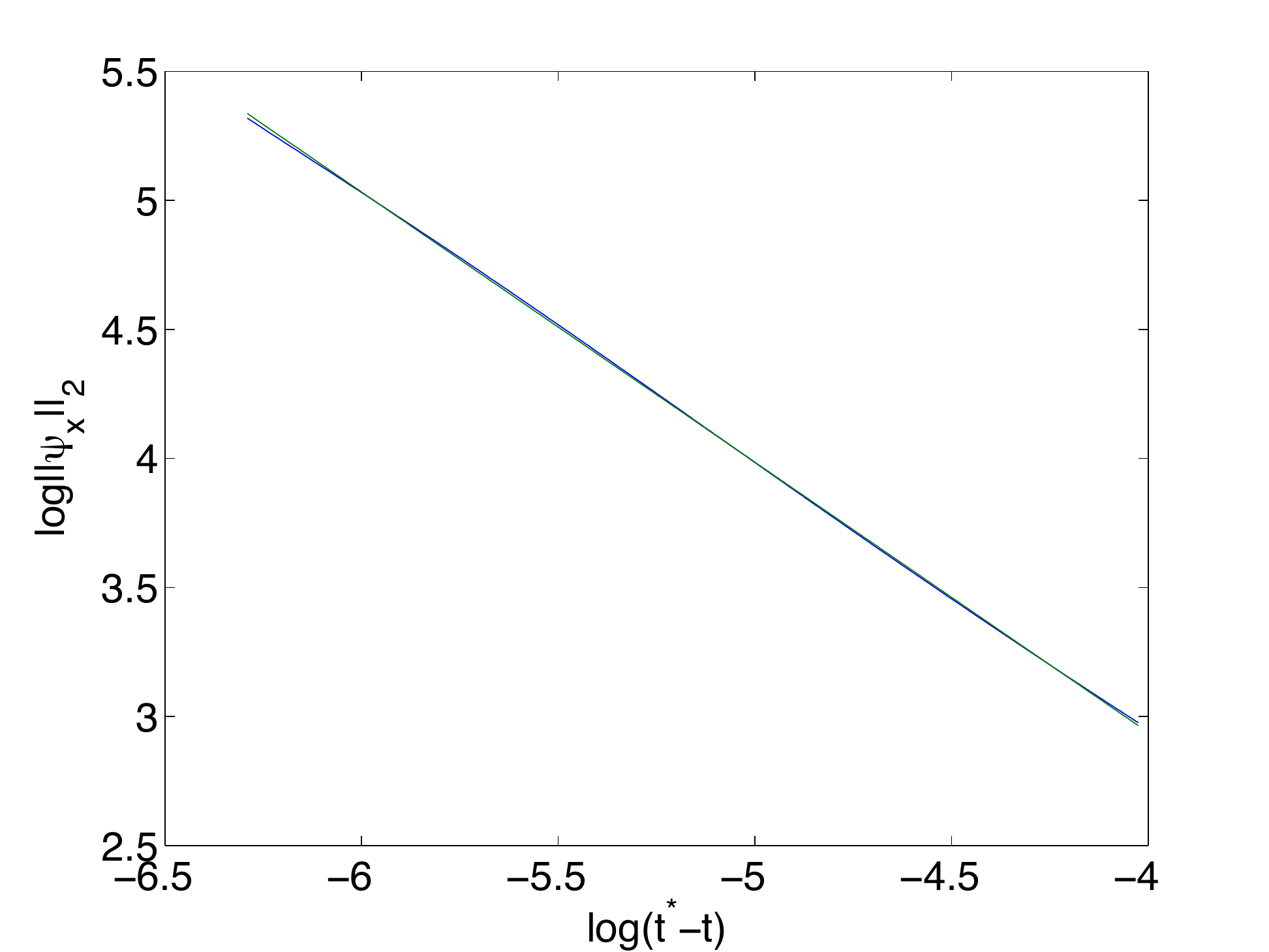}
   \includegraphics[width=0.49\textwidth]{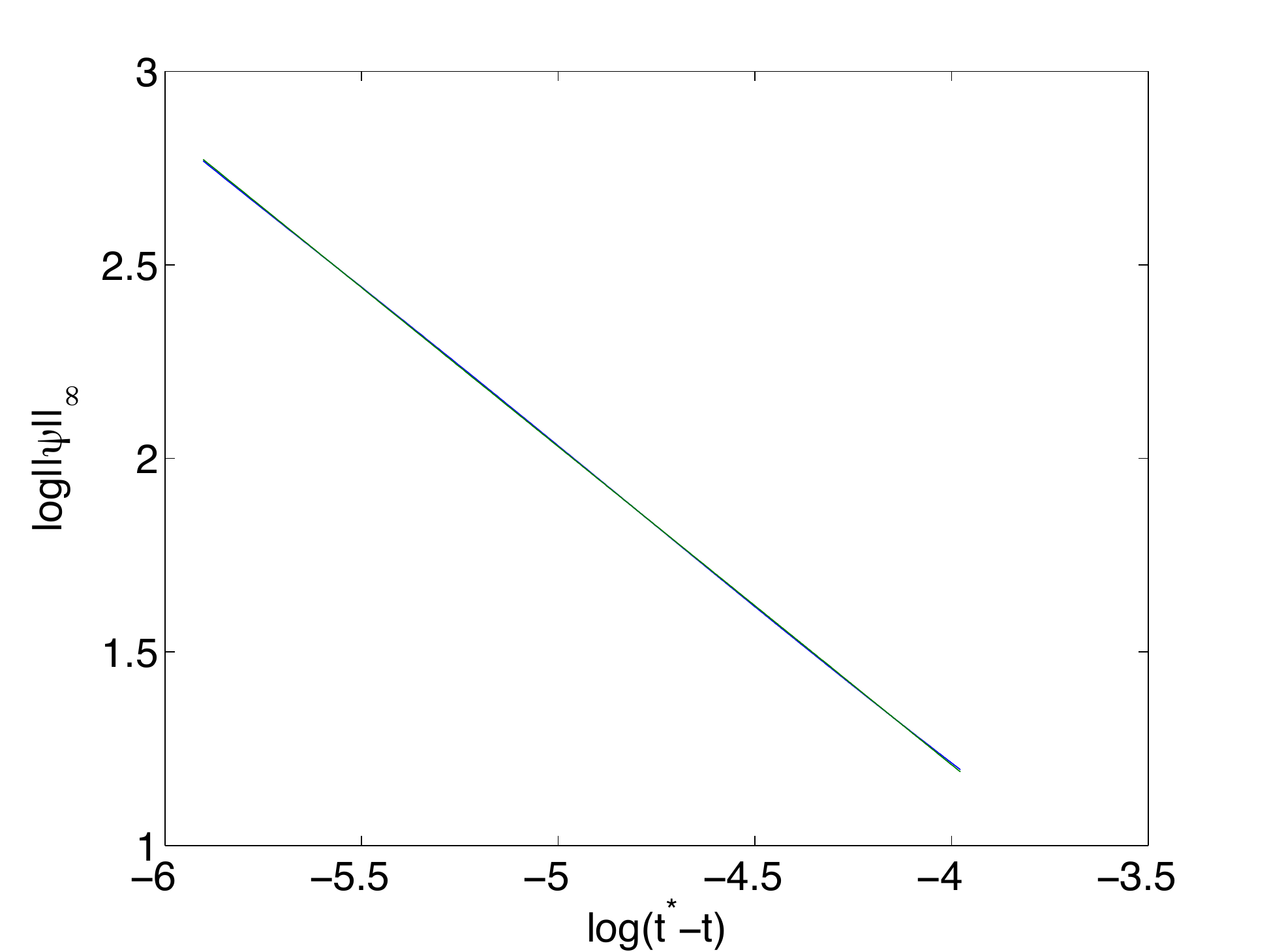}
 \caption{Fitting the logarithms of the $L_{2}$ norm  of $\psi_{x}$ 
 (left) and of the 
 $L^{\infty}$ norm of $\psi$ (right) of the solution to the focusing 
 DS II equation with $\beta=1$ with initial data 
 $\psi_{0}=\exp(-x^{2}-y^{2})$ close to the blow-up.  
The fitted line $\kappa_{1}\ln 
 (t^{*}-t)+\kappa_{2}$ (see the description) is given in green.}
 \label{DSgaussfit}
\end{figure}

The appearance of blow-up is a very subtle and surprising phenomenon in 
DS II type systems since one does not expect blow-up in the cubic hyperbolic NLS equation.
In order to obtain the actual blow-up time, we use 
the optimization algorithm \cite{fminsearch}, which is  
accessible via Matlab as the command \emph{fminsearch}. For $t\approx t^{*}$, we fit for 
the $L_{\infty}$ of the solution and the $L_{2}$ norm of 
$x$-derivative to the expected asymptotic behavior 
(\ref{genscal}). The $L_{\infty}$ norm thereby catches the 
local behavior of the solution close to blow-up, whereas the $L_{2}$ 
norm takes into account the solution on the whole 
computational domain. Thus the consistency of the fitting results 
provides a test of the quality of the numerics. The results of the 
fitting can be seen in Fig.~\ref{DSgaussfit}. 
Fitting $\|\psi_{x} \|_{2}^{2}$ (normalized to 1 at $t=0$) for the last 500 computed time steps 
to $\kappa_{1}\ln (t^{*}-t)+\kappa_{2}$, we 
find $t^{*}=0.2946$, $\kappa_{1}=-1.047$ and $\kappa_{2}=-1.249$. Similarly, 
we get for $\|\psi \|_{\infty}$ the values $t^{*}=.2954$, 
$\kappa_{1}=-.822$ and $\kappa_{2}=-2.078$. Note the 
agreement of the blow-up times which shows the consistency of the 
fitting results within numerical precision. The fitting for the 
$L_{2}$ norm of $\psi_{x}$ agrees very well with the theoretical 
expectation $-1$, whereas the value for the $L_{\infty}$ norm is not 
close to the expected $1/2$. This indicates that we did not get close 
enough to blow-up for lack of resolution to catch the asymptotic 
behavior also locally near the blow-up. Note that these values 
are unchanged within numerical precision if only the last 100 
computed time steps are used for the fitting. 

An interesting question is, whether the logarithmic corrections in (\ref{L2scal}) can also be seen within 
this approach, though certainly not in reliable way due to a lack of 
resolution. 
To test what can be seen with the present code, we do the same fitting as 
above for the last 20 computed time steps since the logarithmic 
corrections will be mainly noticeable for $t\approx t^{*}$, see 
e.g.~the discussion in \cite{KSM}. 
We denote the $L^{2}$ norm of the difference between 
the logarithm of the fitted norm and $\kappa_{1}\ln 
(t^{*}-t)+\kappa_{2}$ as the \emph{fitting error} $\Delta_{2}$. We 
find $\Delta_{2}=0.0061$ for the $L^2$ norm of $\psi_{x}$. If we 
fit the same norms to $\tilde{\kappa}_{1}(\ln 
(t^{*}-t)-\ln\ln|\ln(t^{*}-t)|)+\tilde{\kappa}_{2}$, we get for the 
analogously defined fitting error $\tilde{\Delta}_{2}$ the value 
$0.0035$. Thus there appears to be an indication of logarithmic 
corrections, but this will have to be  checked with higher resolution 
or an adaptive code.


\subsection{Davey-Stewartson solutions in the defocusing case}
\begin{figure}[htb!]
  \includegraphics[width=0.45\textwidth]{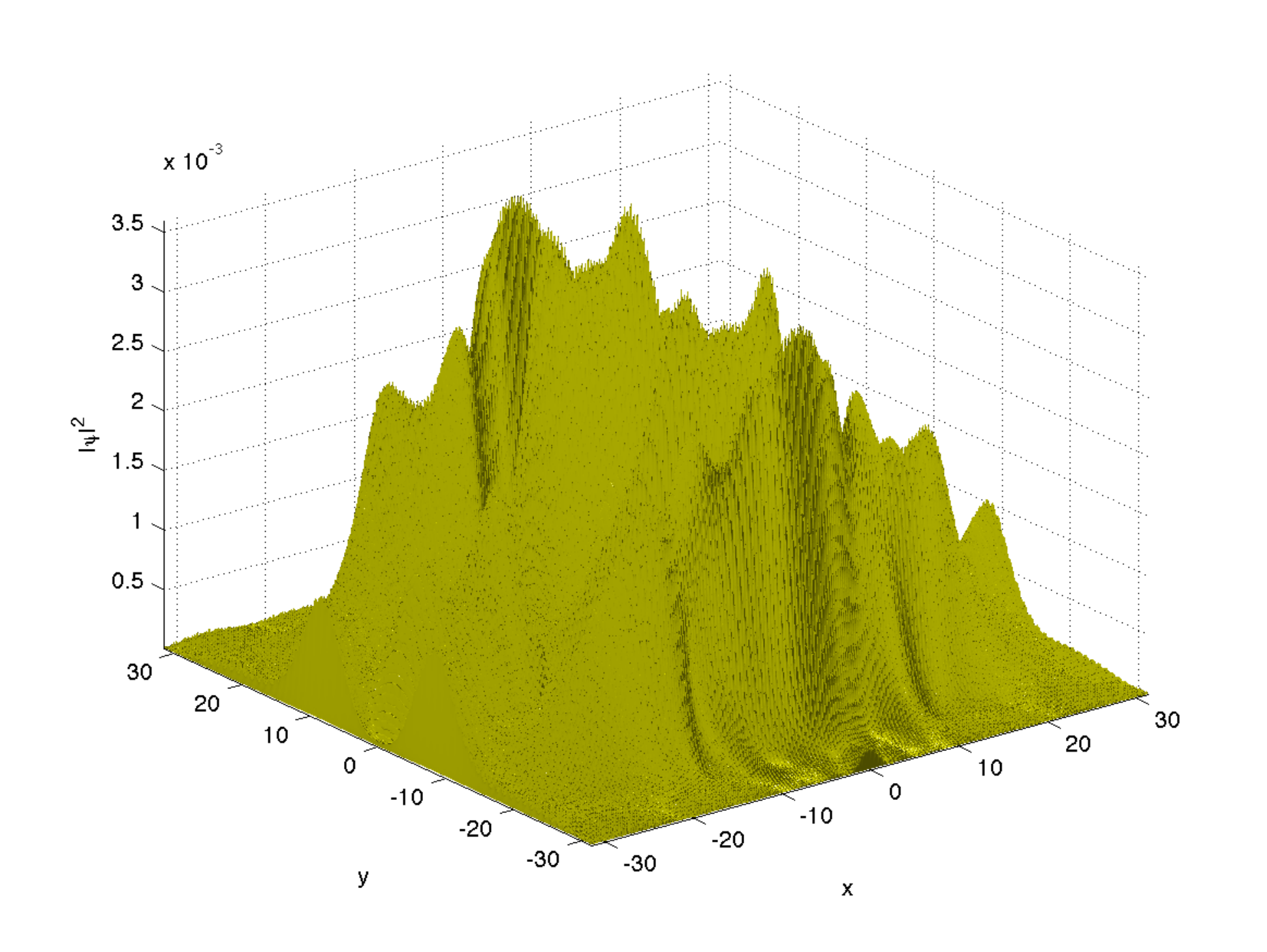}
  \includegraphics[width=0.45\textwidth]{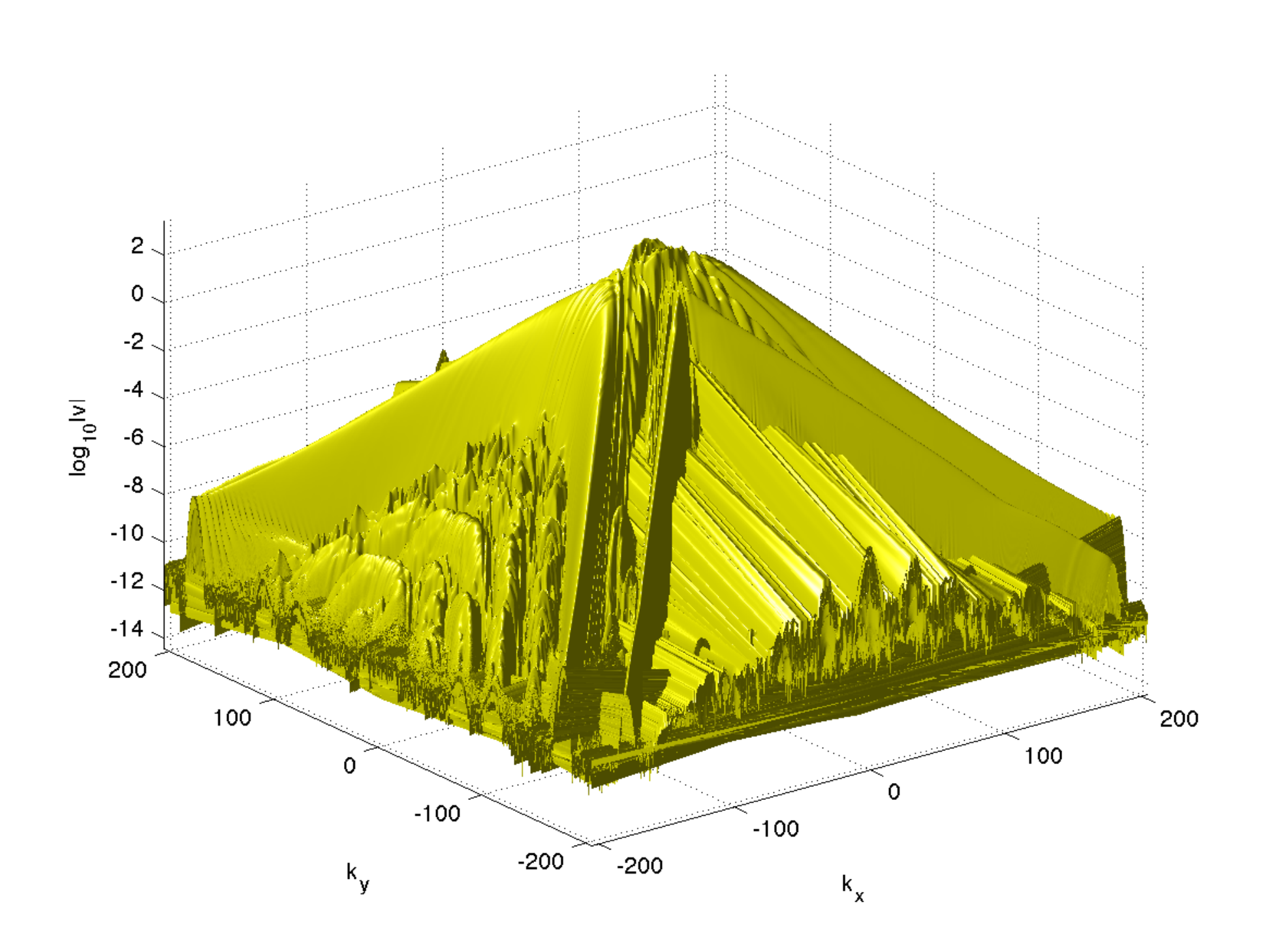}
 \caption{Solution to the hyperbolic NLS equations, i.e., equation 
 (\ref{DSII}) for $\beta=0$ and $\rho=1$ for the initial data 
 $\psi_{0}=\exp(-x^{2}-y^{2})$ at $t=10$ and the corresponding 
 Fourier coefficients.}
 \label{DSdbeta0e01t10}
\end{figure}

In this subsection we will study the same initial data as above for 
the defocusing case, $\rho=1$. For the hyperbolic NLS this implies as 
already discussed an  interchange of $x$ and $y$, i.e., of the 
focusing and of the defocusing direction. If the code is run for 
longer times, the pattern spreads more and more and the $L_{\infty}$ 
norm decreases.  The solution at time $t=10$ can be seen in 
Fig.~\ref{DSdbeta0e01t10}. The computation is carried out for 
$(x,y)\in[-10\pi,10\pi]\times[-10\pi,10\pi]$. Since the initial hump 
reaches the boundaries of the computational domain at large times, 
there are minor effects due to the imposed periodicity conditions 
which are visible in the Fourier coefficients in the same figure. It 
can be recognized, however, that the situation is still numerically 
well resolved. 

\begin{figure}[htb!]
  \includegraphics[width=0.6\textwidth]{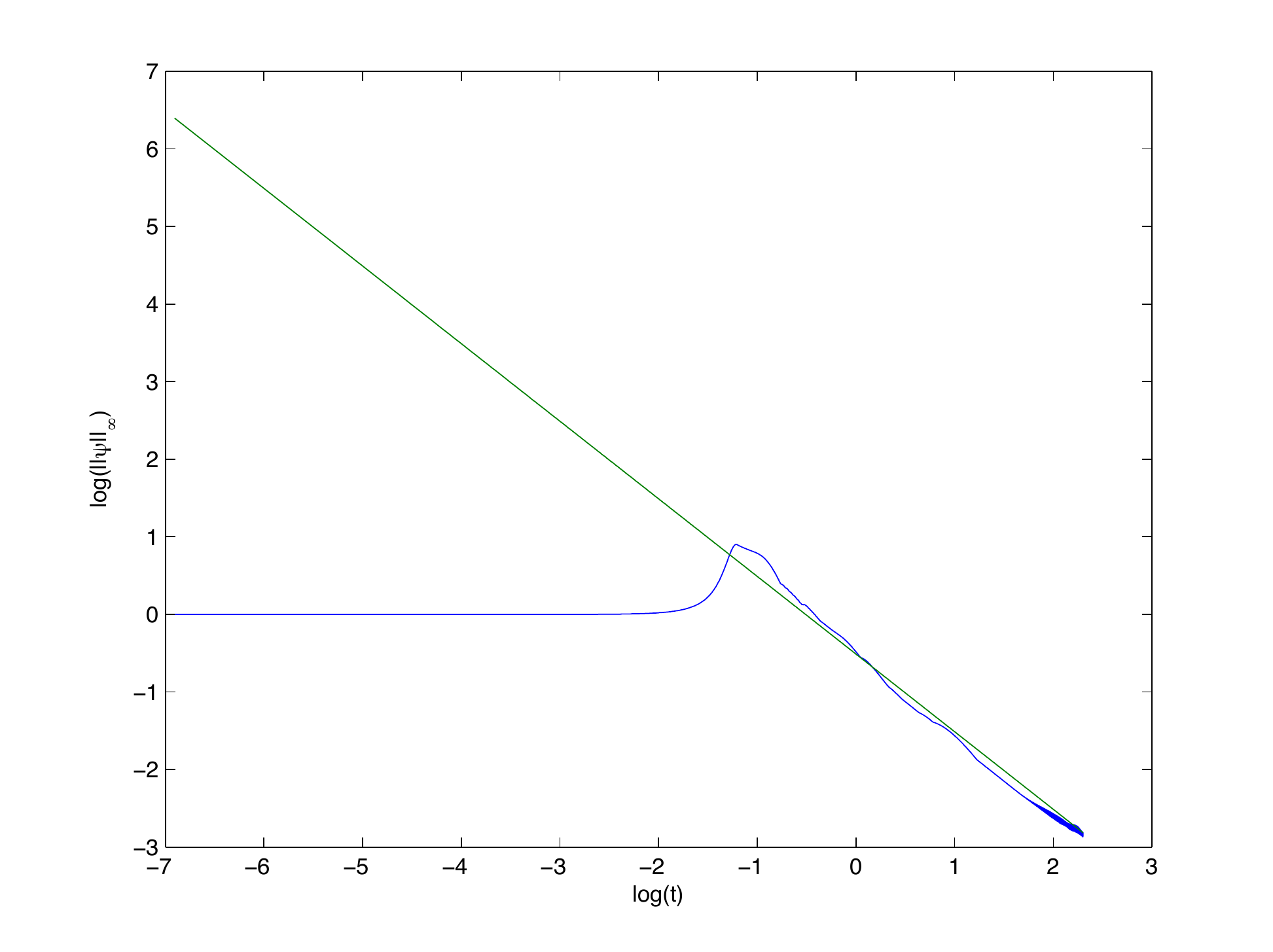}
 \caption{$L_{\infty}$-norm  of the solution to the defocusing DS II equation
 (\ref{DSII}) for $\beta=0$, for the initial 
 data $\psi_{0}=\exp(-x^{2}-y^{2})$ for $\epsilon=0.1$  in dependence 
 of $t$; the green line has slope $-1$ corresponding to a decrease of 
 the $L_{\infty}$ norm proportional to $1/t$.}
 \label{DSdbeta0e01t10max}
\end{figure}

These effects due to the boundary of the computational domain are 
also visible in the $L_{\infty}$ norm in 
Fig.~\ref{DSdbeta0e01t10max} in the form of small oscillations. 
In the logarithmic plot we show also 
a line with slope $-1$ going through the endpoint of the graph of 
$\log ||\psi||_{\infty}$. This suggests that the $L_{\infty}$ norm 
decreases as $1/t$. Since there appear to be no stable solutions to 
DS II (see the conjecture in \cite{MFP} that initial data either lead 
to a blow-up or are radiated away), a $1/t$ decrease of the 
$L_{\infty}$ norm could be observed numerically for all cases without 
blow-up in the limit $t\to\infty$. In the focusing case, this is 
numerically difficult to address which is why we concentrate in this 
context on the defocusing case.

\begin{figure}[htb!]
  \includegraphics[width=\textwidth]{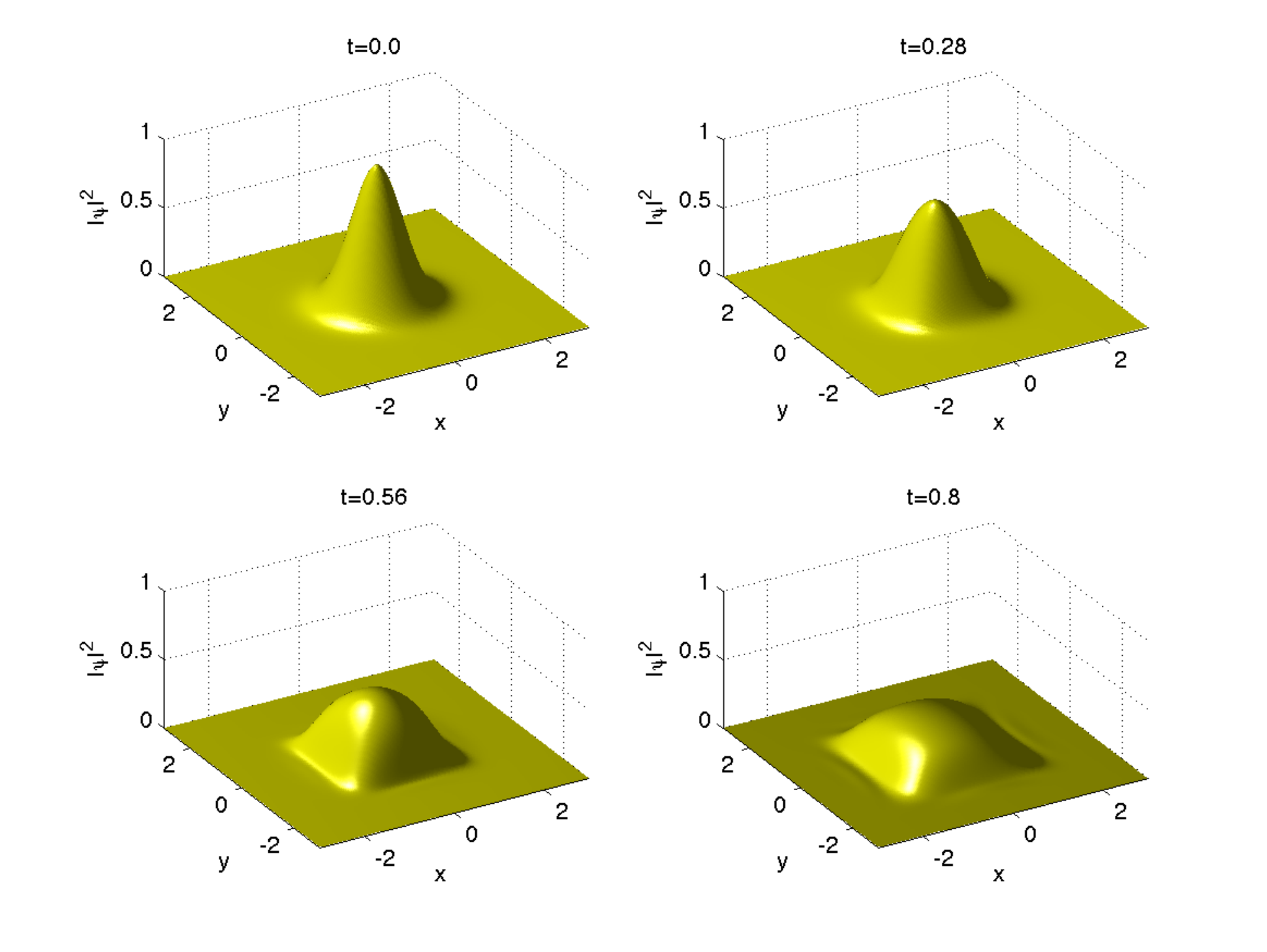}
 \caption{Solution to the defocusing DS II equation
 (\ref{DSII}) for $\beta=0.9$, for the initial 
 data $\psi_{0}=\exp(-x^{2}-y^{2})$ for $\epsilon=0.1$ at different 
 times.}
 \label{DSdbeta09e01t084t}
\end{figure}

For larger values of $\beta$, the behavior of solutions to the 
defocusing DS II equations becomes more similar to what was shown for 
the integrable case in \cite{KR}. For $\beta=0.9$, this can be seen 
for instance in Fig.~\ref{DSdbeta09e01t084t}. The initial pulse is 
defocused to an almost pyramidal shape with a steeping of the 
gradient at the fronts. There are small oscillations in the vicinity 
of the region with strongest gradients. 

\begin{figure}[htb!]
  \includegraphics[width=0.45\textwidth]{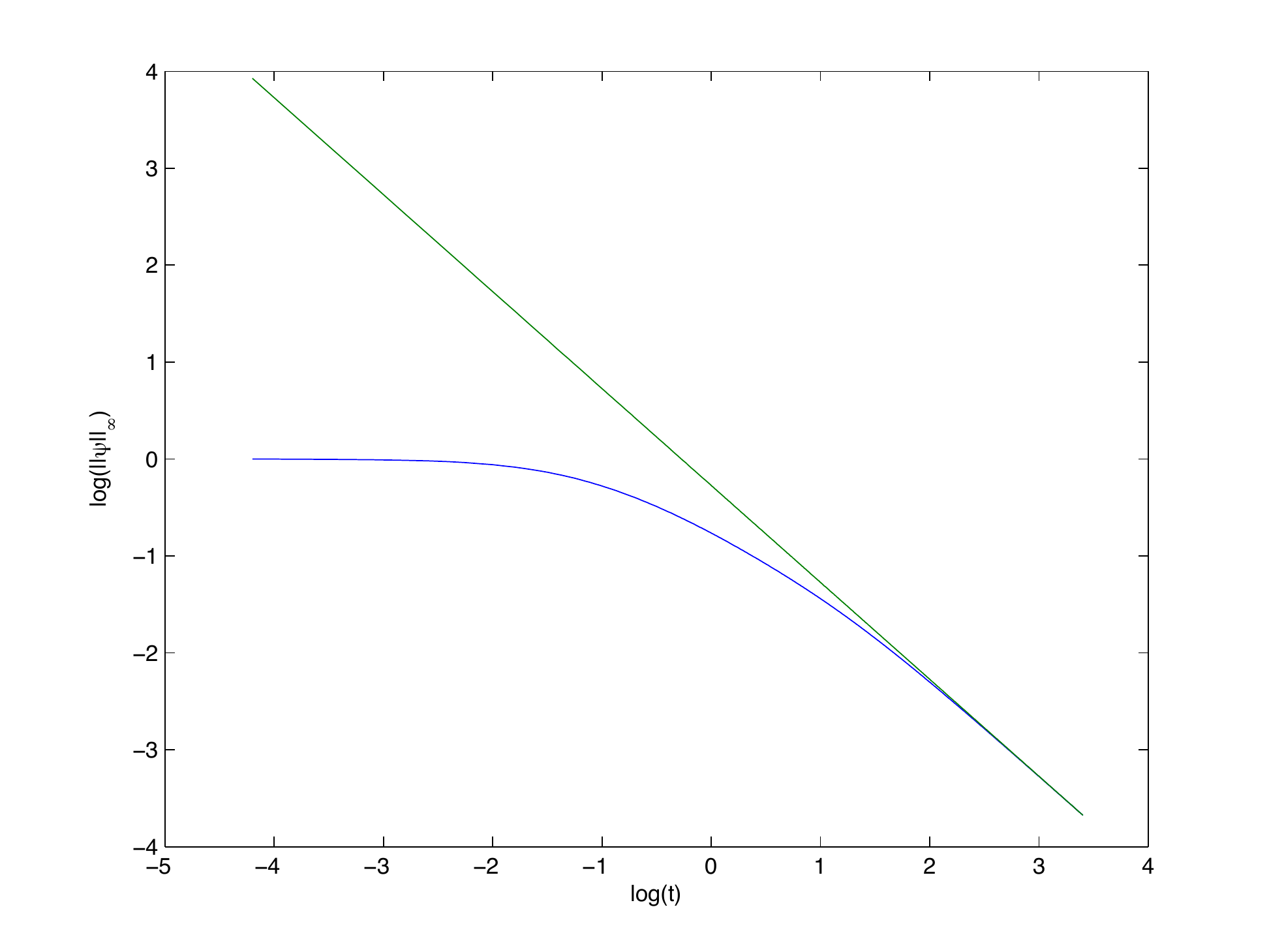}
  \includegraphics[width=0.45\textwidth]{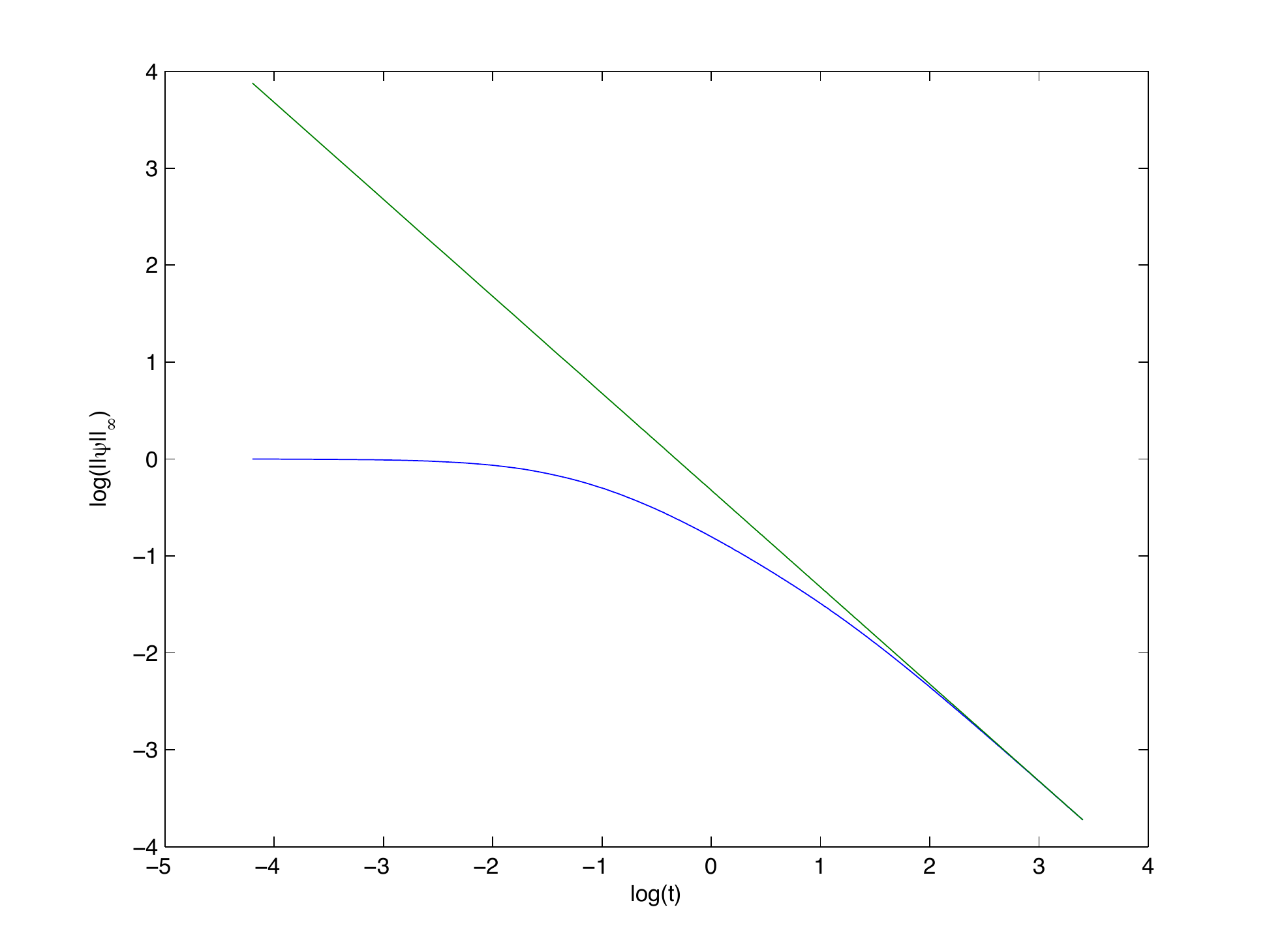}
 \caption{$L_{\infty}$-norm  of the solution to the defocusing DS II equation
 (\ref{DSII}) for $\beta=0.9$ on the left and $\beta=1$ on the right, 
 for the initial 
 data $\psi_{0}=\exp(-x^{2}-y^{2})$ for $\epsilon=0.1$  in dependence 
 of $t$; the green lines have slope $-1$.}
 \label{DSdbeta09e01t30}
\end{figure}

The $L_{\infty}$ norm appears again to decrease as $1/t$ for long 
times as can be seen in Fig.~\ref{DSdbeta09e01t30}. The situation 
is very similar to the integrable case $\beta=1$ in the same figure.

\begin{figure}[htb!]
  \includegraphics[width=\textwidth]{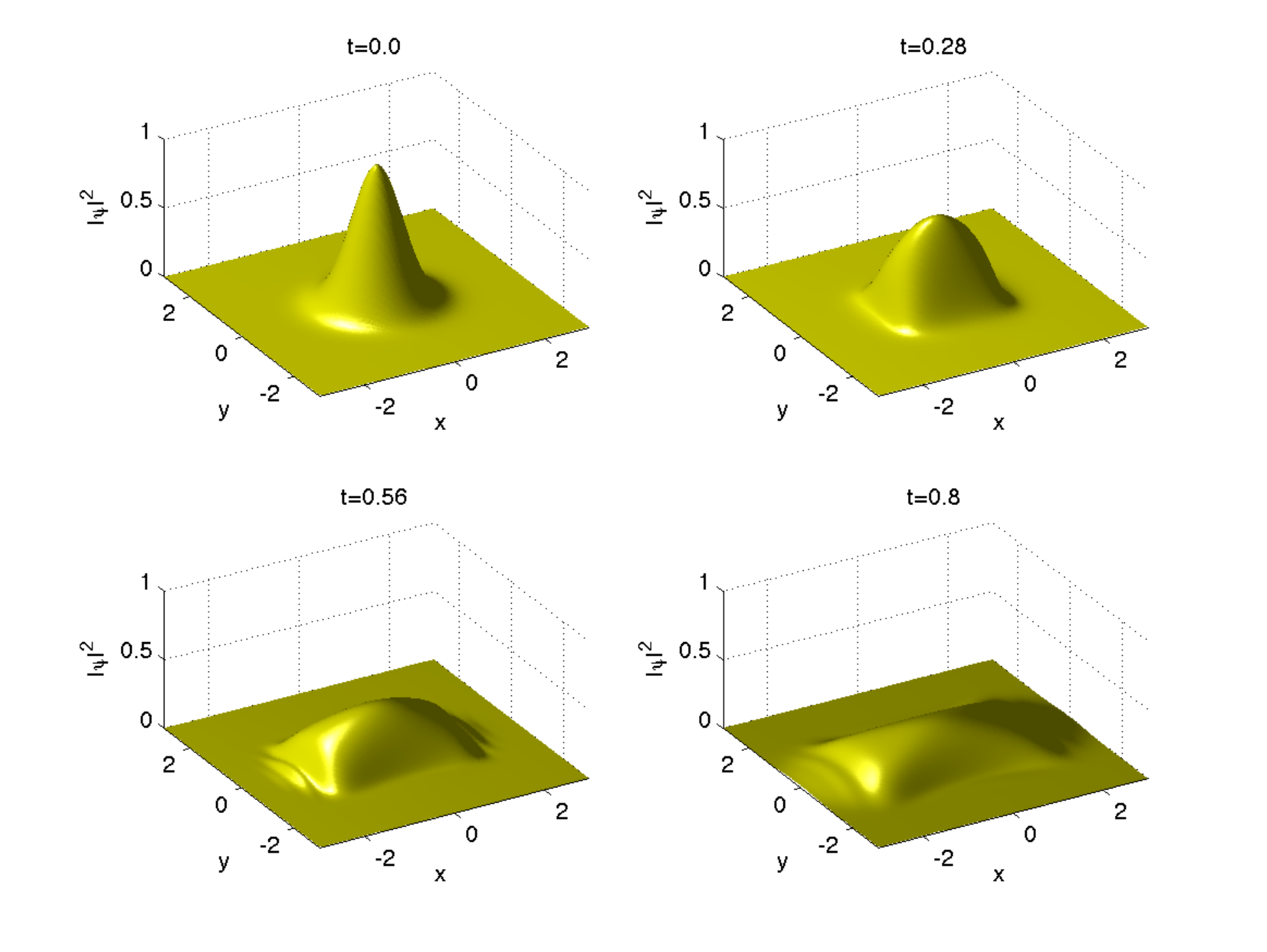}
 \caption{Solution to the defocusing DS II equation
 (\ref{DSII}) for $\beta=2$, for the initial 
 data $\psi_{0}=\exp(-x^{2}-y^{2})$ for $\epsilon=0.1$ at different 
 times.}
 \label{DSdbeta20e01t084t}
\end{figure}

For larger values of $\beta$, the defocusing effect is more 
pronounced in $x$-direction as can be seen in 
Fig.~\ref{DSdbeta20e01t084t}.
The long time behavior is, however, as before, the $L_{\infty}$ norm 
appears to decrease as $1/t$.

\section{Conclusion}
The numerical simulations above together with the results obtained by Inverse Scattering techniques suggest the following conjectures:

1. Cubic hyperbolic NLS: We expect global solutions, with dispersion of the sup norm as $O(1/t).$

2. Defocusing DS II type systems: We expect the same behavior as for the hyperbolic cubic NLS.

3. Focusing DS II type systems. We do not expect finite time blow-up 
except in the integrable case $\beta =1$ for initial data invariant 
under an exchange of the $x$ and $y$. In particular no localized solitary waves should exist when $\beta\neq 1.$

4. Blow-up for the focusing integrable system may occur for initial 
data different to the ones given by Ozawa's construction, say a  
Gaussian with sufficiently high mass and invariant under the 
transformation $x\to y$. Proving such a surprising result should not be easy since the integrable focusing DS system belongs to the "non hyperbolic cubic NLS family" (the dynamic of which is likely to be governed by dispersion) and the standard methods of blow-up for say, the focusing cubic NLS equation (explicit blow-up via a ground state solution and a pseudo-conformal law or virial techniques) clearly do not apply here. Moreover this blow-up, if confirmed, appears to be highly non structurally stable since our simulations suggest  that it does not persist in the non integrable case. In particular, finding a criterion of blow-up such as the one obtained in the $L^2$ critical cubic NLS\footnote{In this case, blow up may occur only for initial data having an $L^2$ norm greater than the ground state one.} is an interesting open question.


\begin{thebibliography}{9}

\bibitem{AC} \textsc{M.J. Ablowitz and P.A. Clarkson}, {\it Solitons, nonlinear evolution equations and inverse scattering}, London Mathematical Society Lecture Notes series {\bf 149}, Cambridge University Press, (1991).

\bibitem{AS} \textsc{M.J. Ablowitz and H. Segur}, {\it On the evolution of packets of water waves}, J. Fluid Mech. {\bf 92} (1979), 691-715.

\bibitem{APP}\textsc{V.A. Arkadiev, A.K. Pogrebkov and M.C. Polivanov}, {\it Inverse scattering transform and soliton solution for Davey-Stewartson II equation},
Physica D {\bf 36} (1089), 188-197.

\bibitem{BR} \textsc{D.J. Benney and G.J. Roskes}, {\it Waves instabilities}, Stud. Appl. Math. {\bf 48} (1969), 377-385.

\bibitem{BB} \textsc{C. Besse and  C.H. Bruneau}, {\it Numerical study of elliptic-hyperbolic Davey-Stewartson system:
dromions simulation and blow-up}, Math. Mod. and Meth. in Appl. Sciences {\bf  8},  (8)
(1998), 1363-1386.

\bibitem{BMS} \textsc{C. Besse, N. Mauser and  H.-P. Stimming}, {\it Numerical study of the Davey-Stewartson system}, M2AN Math. Model. Numer. Anal., {\bf 38},  (6) (2004),  
1035-1054.


\bibitem{CW} \textsc{T. Cazenave, and F.B.  Weissler}, 
{\it Some remarks on the nonlinear Schr\"{o}dinger equation in the critical case}, in Nonlinear semigroups, partial differential equations and attractors (Washington, DC, 1987), 18-29,
Lecture Notes in Math., 1394, Springer, Berlin, 1989. 

\bibitem{Ci}\textsc{R. Cipolatti}, {\it On the existence of standing waves for a Davey-Stewartson system},
Comm. Partial Differential Equations {\bf }17 (1992), no. 5-6, 967-988. 

\bibitem{Ci2}\textsc{R. Cipolatti}, {\it On the instability of ground states for a Davey-Stewartson system}, Ann.Inst. H. Poincar\' e, Phys.Th\' eor. {\bf 58} (1993), 85-104.

\bibitem{Co}\textsc{T. Colin}, {\it Rigorous derivation of the nonlinear Schr\"{o}dinger equation and Davey-Stewartson systems from quadratic hyperbolic systems}, Asymptotic Analysis {\bf 31} (2002), 69-91.

\bibitem{CoLa} \textsc{T. Colin and D. Lannes}, {\it Justification of and long-wave correction to Davey-Stewartson systems from quadratic hyperbolic systems}, Disc. Cont. Dyn. Systems {\bf 11} (1) (2004), 83-100.

\bibitem{DS}\textsc {A. Davey and K. Stewartson}, {\it One three-dimensional packets of water waves}, Proc. Roy. Soc. Lond. A {\bf 338} (1974), 101-110.


\bibitem{DR}\textsc{V.D. Djordjevic and L.G. Redekopp}, {\it On two-dimensional packets of capillary-gravity waves}, J. Fluid Mech. {\bf 79} (1977), 703-714.


    \bibitem{D}\textsc{T. Driscoll}, {\it A composite Runge-Kutta Method for the spectral Solution of semilinear PDEs},
Journal of Computational Physics, {\bf 182} (2002),  357-367.


\bibitem{FPS}\textsc{A. Fokas, D. Pelinovsky and C. Sulem}, {\it Interaction of lumps with a line soliton for the
Davey-Stewartson II equation}, Physica D, {\bf 152-153} (2001), 189-198.



\bibitem{GS} \textsc{J.-M. Ghidaglia and J.-C. Saut}, {\it On the initial value
problem for the 
Davey-Stewartson systems}, Nonlinearity, {\bf 3}, (1990), 475-506. 
\bibitem{GS1} \textsc{J.-M. Ghidaglia and J.-C. Saut}, {\it Non existence of traveling wave solutions to nonelliptic nonlnear Schr\"{o}dinger equations}, J. Nonlinear Sci., {\bf 6} {1996}, 139-145.

\bibitem{GS3} \textsc{J.-M. Ghidaglia and J.-C. Saut}, {\it On the Zakharov-Schulman equations}, in Non- linear Dispersive Waves, L. Debnath Ed., World Scientific, 1992, 83-97.

\bibitem{GS4} \textsc{J.-M. Ghidaglia and J.-C. Saut}, {\it Nonelliptic Schr\"odinger
evolution equations}, 
J. Nonlinear Science {\bf 3}, (1993), 169-195. 


\bibitem{H}\textsc{N.~Hayashi}, {\it Local existence in time of solutions to the elliptic-hyperbolic Davey-Stewartson system without smallness condition on the data}, J. Analyse Math\' ematique {\bf 73}, (1997), 133-164.
\bibitem{NH1}\textsc{ N. Hayashi and H. Hirota}, {\it Local existence in time of small solutions to the elliptic-hyperbolic Davey-Stewartson system in the usual Sobolev space}, Proc. Edinburgh Math. Soc. {\bf 40} (1997),  563-581.

\bibitem{NH2}\textsc{ N. Hayashi and H. Hirota}, {\it Global existence and asymptotic behavior in time of small solutions to the elliptic-hyperbolic Davey-Stewartson system}, Nonlinearity {\bf 9} (1996), 1387-1409. 

\bibitem{Ki}\textsc{O.M. Kiselev}, {\it Asymptotics of solutions of higher -dimensional integrable equations and their perturbations}, J. of Mathematical Sciences, {\bf 138} (6) (2006), 6067-6230.

\bibitem{KP2013}\textsc{C. Klein and R. Peter}, \emph{Numerical study of blow-up in solutions to generalized Korteweg-de Vries equations}. Preprint available at {\tt arXiv:1307.0603}

\bibitem{KSM}\textsc{C.~Klein, C.~Sparber and P.~Markowich}, \textit{Numerical 
study of fractional Nonlinear Schr\"odinger equations}, Preprint available at {\tt arXiv:1404.6262}

\bibitem{KR2013a}\textsc{C. Klein and K. Roidot}, \textit{Numerical study of shock formation in the dispersionless Kadomtsev-Petviashvili equation
and dispersive regularizations}. Phys. D {\bf 265} (2013), 1--25.

\bibitem{KR2}\textsc{C. Klein and K. Roidot},  {\it Numerical Study of the semiclassical limit of the Davey-Stewartson II equations}. Prepint available at {\tt arXiv:1401.4745}.


    \bibitem{KR} \textsc{C.~Klein and K.~Roidot},  \emph{Fourth order time-stepping for Kadomtsev-Petviashvili and 
Davey-Stewartson equations},  SIAM Journal on Scientific Computing Vol. 33, No. 6, DOI: 10.1137/100816663 (2011). 

    \bibitem{KMR} \textsc{ C.~Klein, B.~Muite and K.~Roidot},  \emph{Numerical Study of Blowup in the Davey-Stewartson 
System},  Discr. Cont.  Dyn.  Syst. B, Vol. 18, No. 5, 1361--1387 
(2013).

\bibitem{etna} \textsc{C.~Klein},  \emph{Fourth order time-stepping for low dispersion Korteweg-de 
Vries and nonlinear Schr\"odinger equation},  ETNA Vol. 29 116-135 (2008).

\bibitem{fminsearch}\textsc{J.~C. Lagarias, J.~A. Reeds, M.~H. Wright, and P.~E. Wright,} {\it Convergence properties of the Nelder-Mead simplex method in low dimensions}. 
SIAM J. Optimization {\bf 9} (1998), 112--147.

\bibitem{La} \textsc {D. Lannes},  {\it Water waves: mathematical theory and asymptotics},  Mathematical Surveys and Monographs, vol 188 (2013), AMS, Providence.

\bibitem{Le}\textsc{H. Leblond}, {\it Electromagnetic waves in ferromagnets}, J. Phys. A {\bf 32} (45) (1999), 7907-7932.
 \bibitem{LiPo2} \textsc{F. Linares and G. Ponce}, {\it On the Davey-Stewartson systems}, Ann. Inst. H. Poincar\' e Anal. Non Lin\' eaire, {\bf 10} (1993), 523-548.

\bibitem{MFP}\textsc{M. McConnell, A. Fokas, and B. Pelloni}, Localised coherent Solutions of the DSI and DSII
Equations a numerical Study, Mathematics and Computers in Simulation, 69 (2005), 424-438.

\bibitem{MR} \textsc{N. Mauser and K. Roidot}, {\it Numerical study of the transverse stability of NLS soliton solutions in several classes of NLS type equations}, arXiv:1401.5349v1 [math-ph] 21 Jan 2014.

\bibitem{MeRa} \textsc{F. Merle and P. Rapha\"{e}l}, {\it The blow-up dynamic and upper bound rate for critical nonlinear Schr\"{o}dinger equation}, Ann. of Math (2) {\bf 161} (1) (2005), 157-222.


\bibitem{MRZ}\textsc{S.L. Musher, A.M. Rubenchik and V.E. Zakharov}, {\it Hamiltonian approach to the description of nonlinear plasma phenomena}, Phys. Rep. {\bf 129} (5) (1985), 285-366.

\bibitem{NM}\textsc{A. Newell and J.V. Moloney}, {\it Nonlinear Optics}, Addison-Wesley (1992).

\bibitem{ObRo}\textsc{C. Obrecht and K. Roidot}, In preparation.


\bibitem{Ot1}\textsc{M.Ohta}, {\it Stability and instability of standing waves for the generalized Davey-Stewartson system}, Diff. Int. Eq. {\bf 8} (1995), 1775-1788.

\bibitem{Ot2}\textsc{M.Ohta}, {\it Instability of standing waves for the generalized  Davey-Stewartson system}, Ann. Inst. H. Poincar\' e, Phys. Th\' eor. {\bf 62} (1995), 69-80.

\bibitem{Ot3}\textsc{M.Ohta}, {\it Blow-up solutions and strong instability of standing waves for the generalized  Davey-Stewartson system},  Ann. Inst. H. Poincar\' e, Phys. Th\' eor. {\bf 63} (1995), 111-117.

\bibitem{Oz} \textsc{T. Ozawa}, {\it Exact  blow-up solutions to the Cauchy problem for the Davey-Stewartson systems}, Proc.Roy. Soc. London A, {\bf 436} (1992), 345-349.

\bibitem{PSSW} \textsc{G. Papanicolaou, C. Sulem, P.-L. Sulem and  X.P. Wang}, {\it The focusing singularity of the
Davey-Stewartson equations for gravity-capillary waves}, Physica D {\bf 72} (1994), 61-86.

\bibitem{PS}\textsc{D. Pelinovsky and C. Sulem}, {\it Embedded solitons of the Davey-Stewartson II equation}, CRM
Proceedings and Lecture Notes, Volume {\bf  27} (2001), ed. C. Sulem and I.M.
Sigal (2001), 135-145.

\bibitem{PRKD}\textsc{D.E. Pelinovsky, E.A. Rouvinskaya, O.E. Kurkina and B. Deconincks}, {\it Short-wave transverse instability of line solitons of the two-dimensional hyperbolic nonlinear Schr\"{o}dinger equation}, 
Theoretical and Mathematical Physics, {\bf 179} (1) (2014), 452�461.

\bibitem{Pe}\textsc{P. A. Perry}, {\it Global well-posedness and long time asymptotics for the defocussing Davey-Stewartson II equation in $H^{1,1}(\R^2)$}, 
arXiv:1110.5589v2, 26 Sep 2012.

\bibitem{RT1}\textsc{F. Rousset and N. Tzvetkov}, {Transverse nonlinear instability for some Hamiltonian PDE's}, J. Math.Pures
Appl. {\bf 90} (2008), 550-590.

\bibitem{RT2}\textsc{F. Rousset and N. Tzvetkov}, {\it Transverse nonlinear instability for two-dimensional dispersive models}, Ann. IHP, Analyse Non Lin\' eaire, {\bf 26} (2009), 477-496.


\bibitem{RT3}\textsc{F. Rousset and N. Tzvetkov}, {\it A simple criterion of transverse linear instability for solitary waves}, Math. Res. Lett., {\bf 17} (2010), 157-169

\bibitem{Sc}\textsc{E.I. Schulman}, {\it On the integrability of equations of Davey-Stewartson type}, Theor.Math. Phys. {\bf 56} (1983), 131-136.

\bibitem{SSF}
{\sc C.~Sulem, P.~Sulem, and H.~Frisch}, {\em Tracing complex singularities
  with spectral methods}, J. Comp. Phys., 50 (1983), pp.~138--161.

  \bibitem{SS99} \textsc{C.~Sulem and P.L.~Sulem}, {\em The nonlinear {S}chr{\"o}dinger equation: Self-Focusing and Wave Collapse}. Springer Series in Mathematical Sciences Vol. 139, Springer Verlag 1999.


 \bibitem{Su1}\textsc{L.Y. Sung}, {\it An inverse scattering transform for the Davey-Stewartson equations. I}, J. Math. Anal. Appl. {\bf 183} (1) (1994), 121-154. 
 \bibitem{Su2}\textsc{L.Y. Sung},  {\it An inverse scattering transform for the Davey-Stewartson equations. II}, J. Math. Anal. Appl. {\bf 183} (2) (1994), 289-325. 

 
\bibitem{Su3}\textsc{L.Y. Sung}, {\it An inverse scattering transform for the Davey-Stewartson equations. III}, J. Math. Anal. Appl. {\bf 183} (3) (1994), 477-494. 

\bibitem{Su4}\textsc{L.-Y. Sung}, Long-Time Decay of the Solutions of the Davey-Stewartson II Equations, J. Non-
linear Sci., 5 (1995), pp. 433-452.


\bibitem{Za} \textsc{V.E. Zakharov},
{\it Stability of periodic waves of finite amplitude on the surface
of a deep fluid}, J. Appl. Mech. Tech. Phys. {\bf 2} (1968)
190-194.



\bibitem{ZR}\textsc{V. E. Zakharov and A. M. Rubenchik}, {\it Nonlinear interaction of high-frequency and low frequency
waves}, Prikl. Mat. Techn. Phys., {\bf}  (1972), 84-98.

\bibitem{ZS}\textsc{V. E. Zakharov and, E. I. Schulman}, {\it Degenerate dispersion laws, motion invariants and kinetic equations}, Physica {\bf 1D} (1980), 192-202.

\bibitem{ZS2}\textsc{V. E. Zakharov and, E. I. Schulman}, {\it Integrability of nonlinear systems and perturbation
theory}, in What is integrability? (V.E. Zakharov, ed.), (1991), 185-250, Springer Series
on Nonlinear Dynamics, Springer-Verlag.



\end{thebibliography}
\end{document}